\documentclass[12pt,times,sort&compress,3p]{elsarticle}
\usepackage[utf8x]{inputenc}
\usepackage{amsmath,amssymb}
\usepackage{stackengine}
\usepackage{multirow}
\usepackage{array}
\usepackage[labelformat=simple]{subfig}
\usepackage{float}
\usepackage[skip=1pt,font={normalsize},singlelinecheck=on]{caption}
\usepackage{colortbl}
\usepackage{wasysym}
\usepackage{mathrsfs}
\usepackage{enumerate}

\usepackage{pxfonts}
\usepackage{amsfonts}

\usepackage{lscape}

\usepackage[pdftex,                %
bookmarks         = true,
bookmarksnumbered = true,
pdfpagemode       = None,
pdfstartview      = FitH,
pdfpagelayout     = SinglePage,
colorlinks        = true,
urlcolor          = magenta,
pdfborder         = {0 0 0}
]{hyperref}%

\usepackage{tabularx,ragged2e,booktabs,caption}
\newcolumntype{Y}{>{\centering\arraybackslash}X}

\newtheorem{Pro}{Proposition}

\usepackage{color}
\definecolor{darkred}{rgb}{.7, 0, 0}

\usepackage{geometry}

\journal{\textbf{arXiv}}
\begin{document}
	
	\begin{frontmatter}
		
		\title{Multiple combined gamma kernel estimations for nonnegative data with Bayesian adaptive bandwidths}
	
		\author[rvt,rvt2]{Sobom M. Som\'e\corref{cor1}}
		\ead{sobom.some@uts.bf}		
		\author[rvt1]{C\'elestin C. Kokonendji}
		\ead{celestin.kokonendji@univ-fcomte.fr}
		\author[rvt3]{Smail Adjabi}
			\ead{adjabi@hotmail.com}
		\author[rvt4]{Naushad  A. Mamode Khan}
	\ead{n.mamodekhan@uom.ac.mu}
		\author[rvt5]{Said Beddek}
		\ead{said.beddek@yahoo.fr}	
		\address[rvt]{LST, Universit\'e Thomas SANKARA, Saaba, Burkina Faso}
		\address[rvt2]{LANIBIO, Université Joseph KI-ZERBO,	Ouagadougou, Burkina Faso}
		\address[rvt1]{LMB, Université Bourgogne Franche-Comt\'e, Besan\c{c}on, France}
		\address[rvt4]{Department of Economics and Statistics, University of Mauritius, Reduit, Mauritius}
		\address[rvt3]{Research Unit LaMOS,  University of Bejaia, Bejaia, Algeria}
		\address[rvt5]{Département de Mathématiques, Faculté des Sciences et Sciences Appliquées, Université de Bouira, Bouira, Algeria}
	
		\cortext[cor1]{\textit{Corresponding author:} Laboratoire Sciences et Techniques, Universit\'e Thomas SANKARA, 12 BP 417 Ouagadougou, Burkina Faso.}

\begin{abstract}
A modified gamma kernel should not be automatically preferred to the standard gamma kernel, especially for univariate convex densities with a pole at the origin.  In the multivariate case, multiple combined gamma kernels, defined as a product of univariate standard and modified ones,  are here introduced for nonparametric and semiparametric smoothing of unknown orthant densities with support $[0,\infty)^d$. Asymptotical properties of these multivariate associated kernel estimators are established. Bayesian estimation of adaptive bandwidth vectors using multiple pure combined gamma smoothers, and in semiparametric setup, are exactly derived under the usual quadratic function. The simulation results and four illustrations on real datasets reveal very interesting advantages of the proposed combined approach for nonparametric smoothing, compare to both pure standard and pure  modified gamma kernel versions, and under integrated squared error and average log-likelihood criteria. 	
\end{abstract}
 


\end{frontmatter}

\section{Introduction}

Asymmetric  kernels are known to improve smoothing quality for partially or totally bounded supports; e.g. Scaillet \cite{Sc04}  with inverse and reciprocal inverse Gaussian kernels, \cite{OuimetTolosona2021} using Dirichlet kernels and, \cite{Marchant2013} and \cite{Zoug2018} for uni- and multivariate generalized Birnbaum-Saunders, see also Kakizawa \cite{Ka2021}. However, these kernels induce an additional quantity in the bias that needs reduction via modified versions. See one of the pioneers  Chen~\cite{C99,C00} for beta and gamma kernels, and also \cite{JXK}, \cite{HiruSakudo13,HiruSakudo15}, \cite{Igarashiandkakizawa2014,Igarashiandkakizawa2015}, \cite{MS13} and \cite{Harfouche2020}. In the particular and very useful case of gamma kernels, which of the standard and modified gamma kernels is appropriated for multivariate setup?

In the multivariate setting, let $\mathbf{X}_{1}, \ldots, \mathbf{X}_{n}$ be  independent and identically distributed (iid) $d$-variate random variables with an unknown probability density function (pdf) $f$ on $\mathbb{T}_d^+=[0,\infty)^d$, a subset of $\mathbb{R}^d$ with $d\geq 1$. Then, the multiple standard and modified gamma kernel estimators   $\widetilde{f}_{n}$ and $\widetilde{\widetilde{f}}_{n}$ of $f$ are defined, respectively, for $\mathbf{X}_{i} = (X_{i1},\ldots, X_{id})^{\top}$, $i = 1,\ldots,n$, by
\begin{equation}
\widetilde{f}_{n}(\boldsymbol{x})=\frac{1}{n} \sum_{i=1}^{n} \prod_{j=1}^{d}G_{x_{j},h_{j}}({X}_{ij})\quad \mbox{and}\quad \widetilde{\widetilde{f}}_{n}(\boldsymbol{x})=\frac{1}{n} \sum_{i=1}^{n} \prod_{j=1}^{d}G_{\rho(x_{j};h_{j}),h_{j}}({X}_{ij}),\quad\forall\boldsymbol{x} \in  \mathbb{T}_d^+=[ 0,\infty)^d,\label{gammaestimator}
\end{equation} 
where $\boldsymbol{x}=(x_1,\ldots,x_d)^{\top}$ is the target vector and $\mathbf{h}=(h_1,\ldots,h_d)^{\top}$ is the vector of smoothing parameters with $h_j>0$, $j=1,\ldots,d$; see, e.g. Bouezmarni and Rombouts~\cite{BR10} for multiple modified gamma kernel estimator. Both functions $G_{x,h}(\cdot)$ and $G_{\rho(x;h),h}(\cdot)$ are, respectively, the standard and modified gamma kernels given on the support $\mathbb{S}_{x,h}= [0,\infty)=\mathbb{T}_1$ with  $x\geq 0$ and $h>0$: 
\begin{equation}
G_{x,h}(u)=\frac{u^{x/h}\exp(-u/h)}{\Gamma(1+x/h)h^{1+x/h}}1_{[ 0,\infty)}(u)\quad\mbox{and}\quad G_{\rho(x;h),h}(u)=\frac{u^{\rho(x;h)-1}\exp(-u/h)}{\Gamma(\rho(x;h))h^{\rho(x;h)}}1_{[ 0,\infty)}(u),\label{gam2}
\end{equation} 
where $\Gamma(v)=\int_0^\infty s^{v-1}\exp(-s)ds$ is the classical gamma function with $v>0$,   $1_E$ denotes the indicator function of any given event $E$, and   
\begin{align}
\rho(x;h) = \left\{ \begin{array}{cc}
1+(x/2h)^{2}, & \textrm{if $x\in[0,2h)$}\medskip \\
x/h,  ~~~~~~~~~~~~       & \textrm{if $x\in [2h, \infty)$}.\label{ro}
\end{array} \right.
\end{align} 

The (standard) gamma kernel $G_{x,h}(\cdot)$ appears to be the pdf of the usual gamma distribution, denoted by $\mathcal{G}(1 + x/h,h)$ with shape parameter $1 + x/h$ and scale parameter $h$. The estimators \eqref{gammaestimator} were originally introduced in the univariate case by Chen~\cite{C00} and then used in multivariate case by Bouezmarni and Rombouts \cite{BR10}. Notice that the formulation $(\ref{ro})$ allows  to distinguish the boundary region from the interior one with respect to the bias corrections; see Zhang \cite{zhang2010}, Malec and Schienle \cite{MS13} and Libengué Dobélé-Kpoka and Kokonendji \cite{LC2017} for another choices of the boundary and interior regions and also,  \cite{FK15} for several nonnegative kernels. For instance, Malec and Schienle~\cite{MS13} proposed two refined versions of these modified kernels for further boundary bias reductions and also gave recommendations for the use of all these standard, modified and refined versions. In particular, they showed that the modified gamma kernel should not be systematically preferred  over the standard gamma one, especially for univariate convex densities with pole at $x=0$. One can refer to Kokonendji and Somé \cite{CCKSMS2018} for a general form of the estimators \eqref{gammaestimator}  in multivariate setup. Thus it would be interesting first, to check the nature of each margins and, then, to choose between both standard and modified gamma kernels, accordingly. This compromise kernel denoted {\it multiple combined gamma kernel} is given, for $\boldsymbol{x} \in \mathbb{T}_d^+=[ 0,\infty)^d$ and $\ell=0,1,\ldots,d$, as follows:
\begin{align}
\widetilde{f}_{n}(\boldsymbol{x})&=\frac{1}{n} \sum_{i=1}^{n} \mathbf{G}_{\boldsymbol{x},\mathbf{h};\ell}(\mathbf{X}_{i}) \label{MultiGamma}
\\
&=\frac{1}{n} \sum_{i=1}^{n} \left(\prod_{s=1}^{d-\ell}G_{x_{s},h_{s}}({X}_{is})\right)\left(\prod_{r=1}^{\ell}G_{\rho(x_{r};h_{r}),h_{r}}({X}_{ir})\right).\label{mix_gammaestimator}
\end{align}
We here mention that the combined gamma kernels in \eqref{mix_gammaestimator} have the following (pure) particular cases
\begin{align}
\centering
\left\{ \begin{array}{lc}
\ell=0, ~~~ & \textrm{for multiple (pure) standard,  }\medskip \\
\ell=d,  ~~~     & \textrm{for multiple (pure) modified, }\medskip \\
\ell=1, \ldots, d-1,  ~~~& \textrm{for  multiple (pure) combined.  }\label{parti_cases}
\end{array} \right.
\end{align} 

Alongside with the previous nonparametric approach to estimate unknown pdf $f$, a flexible semiparametric method is also very useful. This approach is a compromise between both pure parametric and  nonparametric methods; see, e.g., Hjort and Glad \cite{HjortGlad1995} and Kokonendji and Somé \cite{KS21} in continuous cases and Kokonendji {\it et al.} \cite{KSB09} and Senga Kiéssé and Mizère \cite{SKM13} in univariate discrete cases. Indeed,  we here assume that any $d$-variate pdf $f$ can be formulated as
\begin{equation}\label{semiP}
f(\boldsymbol{x})=w(\boldsymbol{x};\boldsymbol{\theta})\,p_{d}(\boldsymbol{x};\boldsymbol{\theta}),\;\;\; \forall \boldsymbol{x}\in [0,\infty)^d,
\end{equation}
where $p_d(\cdot;\boldsymbol{\theta})$ is the non-singular parametric part according to a reference $d$-variate distribution with corresponding unknown  parameters $\boldsymbol{\theta}=(\theta_1,\ldots, \theta_d)^\top$ and $w(\cdot;\boldsymbol{\theta}):=f(\cdot)/p_d(\cdot;\boldsymbol{\theta})$ is the unknown  weight function part, to be estimated with the multivariate gamma kernel in \eqref{MultiGamma} or \eqref{mix_gammaestimator}.  However, the choice of a reasonable parametric-start  in \eqref{semiP} is not obvious. See \cite{KP2018} and \cite{CCKTOURE2020} for  recent proposals of relative variability  indexes in count and continuous distributions,  respectively.

The performance of the smoothers in \eqref{mix_gammaestimator} crucially depends on the diagonal bandwidth matrix $\mathbf{h}:=diag(h_1,\ldots, h_d)$ with $d$ real parameters, which is a particular case of the  full symmetric one with $d(d+1)$ real  parameters. The global bandwidth matrices selections such as cross-validation, plug-in and recently global Bayesian are known to have lower performances than their variable  counterparts namely  adaptive and local; see, e.g., \cite{S02} and \cite{zhang2006bayesian}. In nonparametric and semiparametric setups, one can see for example \cite{Zi2015}, \cite{Ercelik2020}, \cite{Some20}, \cite{SK20}, \cite{Belaid2018} and \cite{ZAK14b} for Bayesian adaptive approach using  univariate Birnbaum-Saunders kernel,  scaled inverse Chi-Squared, (univariate and multiple) standard gamma kernel, multiple binomial and  multivariate Gaussian kernel, respectively. See also, \cite{brew2}, \cite{Belaid2016} and \cite{Zi2018} with Bayesian local approaches.

The purpose of this paper is to propose  multiple combined gamma kernels, product of univariate (standard and modified) gamma kernels, for  nonparametric and semiparametric smoothing of unknown pdfs. These relevant estimators are made from univariate gamma kernels \eqref{gam2} chosen according to the shape of each univariate component of  multivariate data. The bandwidth vector shall be investigated from the efficient Bayesian adaptive technique. The inverse gamma is used as prior to  derive the explicit formulas of the posterior distribution and the vector of bandwidths.
The rest of the paper is structured as follows. In Section \ref{sec:est_properties}, we present the combined gamma kernels as multivariate associated kernels, and also give some pointwise properties of the corresponding semiparametric estimator. Section \ref{sec:bayesiennegammamultiple} 
shows the Bayesian adaptive method, where the exact formula of the posterior
distribution and the vector of bandwidth are obtained for the adaptive modified gamma kernel estimator, and in semiparametric setup. In Section \ref{sec:Simulations studies},   simulations studies and comparisons  are conducted  using this Bayesian variable approach with multiple pure standard, pure modified and pure combined gamma kernels for only nonparametric smoothing of unknown pdfs.  Section \ref{sec:Illustration} gives four numerical illustrations on uni-, bi- and trivariate real data sets, including the well-documented Old Faithful geyser. Finally, some concluding remarks are given in Section \ref{sec:conclusion}. Proofs of propositions  are derived in Appendices. Some surface and contour plots of pdf to be estimated and their smoothed densities are illustrated for bivariate case in Supplementary material.

\section{Properties of multiple  combined gamma kernel estimators}
\label{sec:est_properties}

In the following first proposition, we point out the multivariate gamma kernel of \eqref{MultiGamma} as a multivariate associated kernel.

\begin{Pro} \label{Proprod}
	Let $[ 0,\infty)^d = \mathbb{T}_{d}^+$ be the support of the pdf $f$ to be estimated, $\boldsymbol{x} = (x_{1}, \ldots, x_{d})^{\top} \in [ 0,\infty)^d$  the target vector and $\mathbf{h}=(h_1,\ldots,h_d)^{\top}$ the vector of univariate bandwidths with $h_j>0$, $j=1,\ldots,d$. Then, the  multivariate gamma kernel  $\mathbf{G}_{\boldsymbol{x},\mathbf{h};\ell}$ of \eqref{MultiGamma} is a multivariate associated  kernel on support $\mathbb{S}_{\boldsymbol{x}, \mathbf{h}}=[ 0,\infty)^d$ where its random vector has mean vector and covariance matrix  $\boldsymbol{x}+\mathbf{A}\left(\boldsymbol{x}, \mathbf{h}\right)=\left(x_{1} + a_{1}(x_{1}, h_{1}), \ldots, x_{d} + a_{d}(x_{d}, h_{d})\right)^{\top}$ and $ \mathbf{B}\left(\boldsymbol{x}, \mathbf{h}\right)=\mathbf{diag}\left(b_{jj}(x_{j}, h_{j})\right)_{j = 1, \ldots, d}$, respectively. 
In particular,	
\begin{enumerate}[(i)]
\item $\mathbf{A}\left(\boldsymbol{x}, \mathbf{h}\right)=(h_1,\ldots,h_d)$ and  $ \mathbf{B}\left(\boldsymbol{x}, \mathbf{h}\right)=\mathbf{diag}\left(h_j(x_j+h_j)\right)_{j = 1, \ldots, d}$, \textrm{ with $\ell=0$ for multiple standard kernel. }
	\item $\mathbf{A}\left(\boldsymbol{x}, \mathbf{h}\right)$ and $\mathbf{B}\left(\boldsymbol{x}, \mathbf{h}\right)$ are given, for all $j=1,\ldots,d$ by  
 \begin{equation}\label{Sys4}
\begin{aligned}\widetilde{a}_{j}(x_j,h_j)=\dfrac{x^{2}_j+4h_j(h_j-x_j)}{4h_j}1_{[0, 2h_j)}(x_j), \end{aligned} \quad  \widetilde{b}_{j}(x_j,h_j)=\left\{\begin{array}{ll}\begin{aligned}
(x^{2}_j+4h_j^2)/4&   &   &\textrm{if}\: x_j\in[0,2h_j]\\
x_j h_j\qquad &   &   &\textrm{if}\: x_j\in(2h_j,\infty),
\end{aligned}\end{array}\right.
\end{equation}

\textrm{ and $\ell=d$ for multiple modified kernel. }
\item $\mathbf{A}\left(\boldsymbol{x}, \mathbf{h}\right)=\left(h_1,\ldots,h_\ell,\widetilde{a}_{1}(x_1,h_1),\ldots,\widetilde{a}_{d-\ell}(x_{d-\ell},h_{d-\ell})\right)$ and\newline
 $\mathbf{B}\left(\boldsymbol{x}, \mathbf{h}\right)=\mathbf{diag}\left(h_1(x_1+h_1),\ldots,h_\ell(x_\ell+h_\ell),\widetilde{b}_{1}(x_1,h_1),\ldots,\widetilde{b}_{d-\ell}(x_{d-\ell},h_{d-\ell})\right)_{j = 1, \ldots, d}$, with notations in \eqref{Sys4}, and $\ell=1,\ldots,d-1$ for multiple combined kernel.
	\end{enumerate}
\end{Pro} 

\textbf{Proof.} 
The results are obtained directly by calculating the mean vector and covariance matrix of $\mathcal{G}_{\boldsymbol{x},\mathbf{h};\ell}$   the random vector of the pdf $\mathbf{G}_{\boldsymbol{x},\mathbf{h};\ell}$ in  \eqref{MultiGamma} and \eqref{mix_gammaestimator} with \eqref{parti_cases}, for $\ell=0,1,\ldots,d$. $\blacksquare$

Let $\mathbf{X}_1,\ldots,\mathbf{X}_n$ be iid  nonnegative orthant $d$-variate random vectors with unknown pdf $f$ on $\mathbb{T}_{d}^+\subseteq [0,\infty)^d$. The semiparametric estimator of (\ref{semiP}) with (\ref{MultiGamma}) is expressed as follows:
\begin{align}\label{SME}
\widehat{f}_{n}(\boldsymbol{x}) &= p_{d}(\boldsymbol{x};{\boldsymbol{\theta}})\frac{1}{n}\sum_{i=1}^{n}\frac{1}{p_{d}(\mathbf{X}_{i};{\boldsymbol{\theta}})}\mathbf{G}_{\boldsymbol{x},\mathbf{h},\ell}(\mathbf{X}_{i})\nonumber\\
&= \frac{1}{n} \sum_{i=1}^{n}\frac{p_{d}(\boldsymbol{x};{\boldsymbol{\theta}})}{p_{d}(\mathbf{X}_{i};{\boldsymbol{\theta}})}\mathbf{G}_{\boldsymbol{x},\mathbf{h},\ell}(\mathbf{X}_i),\;\;\;\boldsymbol{x}\in\mathbb{T}_d^+, 
\end{align}
where the parameter ${\boldsymbol{\theta}}$  can be known with exact value ${\boldsymbol{\theta}}_0$ or unknown and thus estimated $\widehat{\boldsymbol{\theta}}_n$. Without loss of generality, we give the following results of this section with $\boldsymbol{\theta}= \widehat{\boldsymbol{\theta}}_n$ the known parameter of $\boldsymbol{\theta}$. From (\ref{SME}), we then deduce the nonparametric  associated kernel smoother
\begin{equation}\label{w_tilde}
\widehat{w}_{n}(\boldsymbol{x};\widehat{\boldsymbol{\theta}}_n)= \frac{1}{n}\sum_{i=1}^{n}\frac{1}{p_{d}(\mathbf{X}_{i};\widehat{\boldsymbol{\theta}}_n)}\mathbf{G}_{\boldsymbol{x},\mathbf{h},\ell}(\mathbf{X}_{i})
\end{equation}
of the weight function $x \mapsto w(\boldsymbol{x};{\widehat{\boldsymbol{\theta}}_n})$ which depends  on ${\widehat{\boldsymbol{\theta}}_n}$. 

The following assumptions are required  afterwards for asymptotic properties of the semiparametric estimator \eqref{SME} with multiple combined gamma kernel in \eqref{mix_gammaestimator}.
\begin{description}
	\item\label{a1} ({\bf a1}) The unknown pdf $f$ is in the class $\mathscr{C}^2([0,\infty)^d)$ of twice continuously   differentiable functions. 	
	\item\label{a2} ({\bf a2}) $h_j \rightarrow 0$, $j=1,\ldots,d$ and $n^{-1}\displaystyle \prod_{j=1}^d h_j^{-1/2} \rightarrow 0$ as $n\rightarrow \infty$.	
	\item\label{a3} ({\bf a3}) $\mathbb{I}_{1}=\left \{k \in \{1,\ldots,\ell\}; x_{k} \in [0,2h_{k})\right \}$, $\mathbb{I}_{2}=\left \{m \in \{1,\ldots,d-\ell\}; x_{m} \in [0,2h_{m})\right \}$, and the complementary set $\mathbb{I}^{c}_{}=\left \{j \in \{1,\ldots,d\}~;x_{j} \in [ 2h_{j},\infty)\right \}$ of  $\mathbb{I}_{1}  \times\mathbb{I}_{2}$.
\end{description}

\begin{Pro}\label{PropBiasVarf(x,0)}
	Under Assumption ({\bf a1})  on $f(\cdot)=p_{d}(\cdot;\widehat{\boldsymbol{\theta}}_n)w(\cdot)$, then the multiple combined gamma estimator $\widehat{f}_n (\cdot)=p_{d}(\cdot;\widehat{\boldsymbol{\theta}}_n)\widehat{w}_n(\cdot)$ in (\ref{mix_gammaestimator}) of $f$ satisfies	
	\begin{align*}
\mathrm{Bias}\{\widehat{f}_n(\boldsymbol{x})\}&=
	p_{d}(\boldsymbol{x};\widehat{\boldsymbol{\theta}}_n)\left[w(\boldsymbol{x})-f(\boldsymbol{x})\{p_{d}(\boldsymbol{x};\widehat{\boldsymbol{\theta}}_n)\}^{-1}+\displaystyle \sum_{r=1}^{d-\ell}h_r \frac{\partial w}{\partial x_r}(\boldsymbol{x})+\displaystyle \sum_{r=1}^{d-\ell}\frac{1}{2}\left(x_r h_r+2h_r^2 \right)\frac{\partial^2 w}{\partial x_r^2}(\boldsymbol{x})\right.\\
&\quad \left.+\frac{1}{2}\displaystyle \sum_{s=1}^{\ell}x_s h_s \frac{\partial^2 w}{\partial  x_s^2}(\boldsymbol{x})\right] +\left(1+o\left\{\displaystyle \sum_{j=1}^{d}h_j^2 \right\}\right), 
\end{align*}  
for any $\boldsymbol{x}\in[0,\infty)^d$.  Moreover, if  ({\bf a2}) and ({\bf a3}) hold, then there exists $\lambda_{k}>0$, for $k\in \mathbb{I}_1$  such that
\begin{align}\label{Varf}
	\mathrm{var}\{\widehat{f}_n(\boldsymbol{x})\}&=\frac{1}{n}f(\boldsymbol{x})[p_{d}(\boldsymbol{x};\widehat{\boldsymbol{\theta}}_n)]^{-2}\displaystyle \prod_{k \in \mathbb{I}^{}_{1}}\left(\frac{\Gamma(2\lambda_{k}+1)}{2^{2\lambda_{k}+1}\Gamma(\lambda_{k}+1)}h_{k}^{-1}\right)\prod_{s \in \mathbb{I}^{}_{2}}\left(\frac{\Gamma(2\lambda_{s}^2+1)}{2^{2\lambda_{s}^2+1}\Gamma(\lambda_{s}^2+1)}h_{s}^{-1}\right)\nonumber\\
	&\quad \times \prod_{j \in \mathbb{I}_{}^c}\left(\frac{1}{2\pi^{1/2}}h_j^{-1/2}x_j^{-1/2}\right)+o\left(n^{-1}\displaystyle\prod_{j=1}^{d}h_j^{-1/2}\right).
\end{align}
\end{Pro}

\section{Bayesian adaptive bandwidths for  multiple combined gamma kernel estimators}
\label{sec:bayesiennegammamultiple}

In this section, we first provide the Bayesian adaptive approach to select the variable smoothing parameter, suitable  for the modified gamma kernel estimator \eqref{SME}  of unknown pdf, in semiparametric context.
Following the multiple standard gamma case of \cite{SK20} and \cite{KS21}, this modified gamma kernel estimator is constructed from \eqref{SME} with \eqref{gam2} and \eqref{ro} by considering a variable bandwidth vector $\mathbf{h}_{i}=(h_{i1},\ldots,h_{id})^\top$ for each observation $\mathbf{X}_{i}=(X_{i1},\ldots,X_{id})^\top$ instead of  the fixed bandwidth vector $\mathbf{h}=(h_{1},\ldots,h_{d})^{\top}$. Thus, the smoothing parameter $\mathbf{h}_{i}$ is considered as a random vector with a prior distribution.


The semiparametric estimator \eqref{SME} with parametrization \eqref{ro} for multiple modified gamma kernel  and variable vector of bandwidth $\mathbf{h}_{i}$ is written in $\boldsymbol{x}=(x_{1},\ldots,x_{d})^\top\in [0,\infty)^d$ as 
\begin{align}
\widehat{f}_{n}(\boldsymbol{x})&=p_{d}(\boldsymbol{x};\widehat{\boldsymbol{\theta}}_n)\frac{1}{n}\sum_{i=1}^{n}\frac{1}{p_{d}(\mathbf{X}_{i};\widehat{\boldsymbol{\theta}}_n)} \prod_{\ell=1}^{d}G_{\rho(x_{\ell},h_{i\ell}),h_{i\ell}}(X_{i\ell})\nonumber\\
&= \frac{1}{n} \sum_{i=1}^{n}\frac{p_{d}(\boldsymbol{x};\widehat{\boldsymbol{\theta}}_n)}{p_{d}(\mathbf{X}_{i};\widehat{\boldsymbol{\theta}}_n)}\prod_{\ell=1}^{d}G_{\rho(x_{\ell},h_{i\ell}),h_{i\ell}}(X_{i\ell}), \label{variableestimator}
\end{align}
where $\widehat{\boldsymbol{\theta}}_n$ is the estimated parameter of $\boldsymbol{\theta}$.  

Then, the leave-one-out kernel estimator of $f(\mathbf{X}_{i})$  deduced from ($\ref{variableestimator}$) is
\begin{equation}
\widehat{f}_{n,\mathbf{h}_i,-i}(\mathbf{X}_i):=\frac{p_{d}(\mathbf{X}_i;\widehat{\boldsymbol{\theta}}_n)}{n-1}\sum_{\ell=1,\ell\neq i}^{n}
\frac{1}{p_{d}(\mathbf{X}_{\ell};\widehat{\boldsymbol{\theta}}_n)} \prod_{\ell=1}^{d}G_{\rho(x_{\ell},h_{i\ell}),h_{i\ell}}(X_{j\ell}),\label{equ5}
\end{equation} 
where  $\{\mathbf{X}_{-i}\}$ stands for the set of observations excluding $\mathbf{X}_{i}$. Let $\pi(\mathbf{h}_{i})$ be the prior distribution of $\mathbf{h}_{i}$, then the posterior distribution for each variable bandwidth vector $\mathbf{h}_{i}$ given $\mathbf{X}_{i}$ provided from the Bayesian rule is expressed as follows
\begin{equation}
\pi(\mathbf{h}_{i}\mid\mathbf{X}_{i})=\frac{\widehat{f}_{n,\mathbf{h}_i,-i}(\mathbf{X}_i) \pi (\mathbf{h}_{i})}{\int_{\chi}\widehat{f}_{n,\mathbf{h}_i,-i}(\mathbf{X}_i) \pi (\mathbf{h}_{i})d\mathbf{h}_{i}},\label{bayesrule}
\end{equation}
where  $\chi$ is the space of positive vectors. The Bayesian estimator  $\widehat{\mathbf{h}}_{i}$ of $\mathbf{h}_{i}$ is obtained through the usual quadratic loss function as
\begin{equation}
\widehat{\mathbf{\mathbf{h}}}_{i}=\mathbb{E}\left(\mathbf{h}_{i}\mid\mathbf{X}_{i})=(\mathbb{E}(h_{i1}\mid\mathbf{X}_{i}),\ldots,\mathbb{E}(h_{id}\mid\mathbf{X}_{i})\right)^{\top}. \label{equ6}
\end{equation}

We assume that each component $h_{i \ell}=h_{i \ell}(n)$, $\ell=1,\ldots,d$, of $\mathbf{h}_{i}$ has the univariate inverse gamma prior  $\mathcal{IG}(\alpha,\beta_{\ell})$ distribution with same shape parameters $\alpha > 0$ and, eventually, different scale parameters  $\beta_{\ell} > 0$ such that $\boldsymbol{\beta}=(\beta_1,\ldots,\beta_d)^\top$. We recall that the pdf of $\mathcal{IG}(\alpha,\beta_\ell)$ with $\alpha,\beta_\ell>0$ is defined by 
\begin{equation}
IG_{\alpha,\beta_\ell}(u)=\frac{\beta_\ell^{\alpha}}{\Gamma(\alpha)} u^{-\alpha-1} \exp(-\beta_\ell/u)1_{(0, \infty)}(u), \;\;\ell=1,\ldots,d,\label{prior}
\end{equation} 
where $\Gamma(\cdot)$ is the usual gamma function. From those considerations, the closed form of the posterior density and the Bayesian estimator of the vector $\mathbf{h}_{i}$ are given in the following proposition.

\begin{Pro}\label{theo1} 
	For fixed $i \in \{1,2,\ldots,n\}$, consider each observation $\mathbf{X}_{i}=(X_{i1},\ldots,X_{id})^{\top}$ with its corresponding vector  $\mathbf{h}_{i}=(h_{i1},\ldots,h_{id})^{\top}$ of univariate bandwidths and defining the subset  $\mathbb{I}_{\mathbf{X}_i}=\left \{k \in \{1,\ldots,d\}; X_{ik} \in [0,2h_{ik})\right \}$ and $\mathbb{I}^{c}_{\mathbf{X}_i}=\left \{\ell \in \{1,\ldots,d\}~;X_{i\ell} \in [ 2h_{i\ell},\infty)\right \}$ its complementary set. Using the inverse gamma prior $IG_{\alpha,\beta_{\ell}}$ of (\ref{prior}) for  each component $h_{i\ell}$ of $\mathbf{h}_{i}$ in the multiple gamma estimator (\ref{variableestimator}) with $\alpha>1/2$ and $\boldsymbol{\beta}=(\beta_1,\ldots,\beta_d)^\top\in(0,\infty)^d$, then:
	
	(i) there exists $\lambda_{ik}>0$ for $k\in \mathbb{I}_{\mathbf{X}_i}$ such that the posterior density (\ref{bayesrule}) is the following weighted sum of inverse gamma 
	\begin{align*}
	\pi(\mathbf{h}_{i}\mid\mathbf{X}_{i})&=\frac{p_{d}(\mathbf{X}_i;\widehat{\boldsymbol{\theta}}_n)}{D_{i}(\alpha,\boldsymbol{\beta})}\sum_{j=1,j\neq i}^{n}\frac{1}{p_{d}(\mathbf{X}_{j};\widehat{\boldsymbol{\theta}}_n)}\left(\prod_{k \in \mathbb{I}_{\mathbf{X}_i}}A_{ijk}(\alpha,\beta_k)\,IG_{\lambda_{ik}+\alpha+1,X_{jk}+\beta_{k}}(h_{ik})\right) \\\nonumber
	&\quad \times \left(\prod_{\ell \in \mathbb{I}^{c}_{\mathbf{X}_i}} B_{ij\ell}(\alpha,\beta_\ell)\,IG_{\alpha+1/2,C_{ij\ell}(\beta_\ell)}(h_{i\ell})\right),
	\end{align*}
	with $A_{ijk}(\alpha,\beta_k)= [ \Gamma (\lambda_{ik}+ \alpha +1)X_{jk}^{\lambda_{ik}} ]/[\beta_{k}^{-\alpha} \Gamma(\lambda_{ik}+1 )(X_{jk}+\beta_{k})^{\lambda_{ik}+\alpha+1}]$, $B_{ij\ell}(\alpha,\beta_\ell)= [X^{-1}_{j\ell}\Gamma(\alpha +1/2)]/(\beta_{\ell}^{-\alpha}X_{i\ell}^{-1/2}\sqrt{2\pi}[C_{ij\ell}(\beta_\ell)]^{\alpha +1/2})$, $C_{ij\ell}(\beta_\ell)= X_{i\ell}\log(X_{i\ell}/X_{j\ell})+X_{j\ell}-X_{i\ell}+\beta_{\ell}$,  and $D_{i}(\alpha,\boldsymbol{\beta})=p_{d}(\mathbf{X}_i;\widehat{\boldsymbol{\theta}}_n)\sum_{j=1,j\neq i}^{n}\left(p_{d}(\mathbf{X}_{j};\widehat{\boldsymbol{\theta}}_n)\right)^{-1}\left(\prod_{k \in \mathbb{I}_{\mathbf{X}_i}}A_{ijk}(\alpha,\beta_k)\right)\left(\prod_{\ell \in \mathbb{I}^{c}_{\mathbf{X}_i}} B_{ij\ell}(\alpha,\beta_\ell)\right)$;
	
	(ii) under the quadratic loss function, the Bayesian estimator $\widehat{\mathbf{\mathbf{h}}}_{i}=\left(~\widehat{h}_{i1},\ldots,\widehat{h}_{id}\right)^{\top}$  of $\mathbf{h}_{i}$, introduced in (\ref{equ6}), is 
	\begin{align*}
	\widehat{h}_{im} &= \frac{p_{d}(\mathbf{X}_i;\widehat{\boldsymbol{\theta}}_n)}{D_{i}(\alpha,\boldsymbol{\beta})}\sum_{j=1,j\neq i}^{n}\frac{1}{p_{d}(\mathbf{X}_{j};\widehat{\boldsymbol{\theta}}_n)}\left(\prod_{k \in \mathbb{I}_{\mathbf{X}_i}}A_{ijk}(\alpha,\beta_k)\right)
	\left(\prod_{\ell \in \mathbb{I}^{c}_{\mathbf{X}_i}} B_{ij\ell}(\alpha,\beta_\ell)\right) \\\nonumber
	&\quad \times \left(\frac{X_{jm}+\beta_{m}}{\lambda_{ik}+\alpha}1_{[0,2h_{im})}(X_{im}) + \frac{C_{ijm}(\beta_m)}{\alpha-1/2}1_{[2h_{im},\infty)}(X_{im})\right),
	\end{align*}
	for $m=1,2,\ldots,d,$ with the previous notations of   $A_{ijk}(\alpha,\beta_k)$, $B_{ij\ell}(\alpha,\beta_\ell)$, $C_{ijm}(\beta_m)$ and $D_{i}(\alpha,\boldsymbol{\beta})$.
\end{Pro}

One can notice that $\lambda_{ik} \rightarrow 0$ provides the Bayesian adaptive bandwidth vector for the nonparametric multiple gamma kernel smoother of \cite{SK20}. Furthermore,  $p_{d}(\cdot;{\boldsymbol{\theta}}) \equiv 1$ gives the nonparametric setup. Thus, the Bayesian adaptive selector for the new multiple combined gamma kernel can be easily deduced as combination of the standard and modified cases. Similarly to \cite{SK20} and \cite{KS21} for nonparametric and semiparametric approaches respectively,   we have to consider  $\alpha = \alpha_n = n^{2/5}>2$ and $\beta_\ell=1>0$, $\ell=1,\ldots,d$ to ensure the convergence of the variable bandwidths to zero and also in practice. Although, these previous choices do not  give necessarily  the best smoothing quality.

\section{Simulation studies}
\label{sec:Simulations studies}

In this section, all numerical studies are performed only in nonparametric context with  multiple standard, modified and combined gamma kernel smoothers \eqref{gammaestimator}, \eqref{gam2} and \eqref{ro}. Thus, we provide simulation results, conducted for evaluating the performance of the proposed approach, namely  Bayesian adaptive bandwidth selection for these three types of multiple  gamma kernel estimators  \eqref{mix_gammaestimator} with \eqref{parti_cases}.  Computations have been run on PC 2.30 GHz by using the \textsf{R} software \cite{R20}. 
The following numerical studies have two objectives with respect to the standard gamma kernel method of \cite{{SK20},{KS21}}. Firstly, we investigate the ability of our proposed method (\ref{variableestimator}) with pure modified and combined gamma kernels \eqref{mix_gammaestimator} with \eqref{parti_cases} to produce good nonparametric estimate  of unknown true densities on $[0,\infty)^d$. Finally, we compare execution times needed for computing all Bayesian  procedure with smoothers \eqref{gammaestimator} and  \eqref{mix_gammaestimator}. 

In fact, the implementation of the Bayesian adaptive  approach  with inverse gamma priors requires the choice of parameters $\alpha$ and $\boldsymbol{\beta}=(\beta_1,\ldots,\beta_d)^\top$. Throughout this section, we take $\alpha=n^{2/5}>2$ and $\beta_\ell = 1$ for all $\ell=1,2,\ldots,d$. 
The efficiency of the smoothers shall be examined through the empirical estimate $\widehat{ISE}$ of the integrated squared errors (ISE): 
$$ISE:=\displaystyle{\int}_{\left[0,\infty\right)^d}\left\lbrace\widetilde{f}(\boldsymbol{x})- f(\boldsymbol{x})\right\rbrace^{2} d\boldsymbol{x},$$ 
where $\widetilde{f}$ is   the variable multivariate associated gamma kernel estimator $\widehat{f}_n$ of (\ref{variableestimator}) for our adaptive Bayesian method. All the $\widehat{ISE}$ values are here computed with the number of replications $N=100$.

The boundary region in \eqref{ro} are taken small enough (in practice $h_\ell\approx 10^{-8}$ for all $\ell=1,\ldots,d$)  to contain very little or no observations, and with respect to the modified  and combined gamma kernels \eqref{mix_gammaestimator}. 

\subsection{Univariate case study}\label{SS_univ}

We here consider four scenarios denoted A, B, C and D  to simulate nonnegative datasets  with respect to the support of both standard or modified gamma kernel (i.e. $\mathbb{S}_{x,h}= [0,\infty)=\mathbb{T}$). These scenarios have also been chosen to compare the performances of both smoothers \eqref{gammaestimator} and \eqref{variableestimator} on dealing with  convex, unimodal or multimodal  distributions.   
	\begin{itemize}		
		\item Scenario A is generated  using  the gamma distribution  
		$$f_{A}(x)=2  \exp(-2 x)1_{[0, \infty)}(x);$$		
		\item  Scenario B comes from   a mixture of three gamma densities 
		$$f_{B}(x)=\left(\frac{3}{10}  \frac{x \exp (-3 x)}{3^{-2}\Gamma\left(2\right)}+\frac{2}{5}  \frac{x^{7} \exp (-3 x)}{3^{-8}\Gamma\left(8\right)}+\frac{3}{10}  \frac{x^{10} \exp (-3 x)}{3^{-11}\Gamma\left(11\right)}\right)1_{[0, \infty)}(x);$$		
		\item Scenario  C is the Weibull density 
		$$f_C(x)=\frac{3}{2}\left(\frac{x}{2}\right)^{2} \exp\left(-[x/ 2]^{3}\right)1_{[0, \infty)}(x);$$		
		\item Scenario D is from a mixture of three Erlang distributions 
		$$f_D(x)=\left(\frac{1}{3}\frac{ x \exp(-2.5 x)}{2.5^{-2}}+\frac{1}{3}\frac{2.5^{8} x^{7} \exp(-2.5 x)}{7 !}+\frac{1}{3}\frac{2.5^{25} x^{24} \exp(-2.5 x)}{24 !}\right)1_{[0, \infty)}(x).$$ 		
	\end{itemize}

Table~\ref{computtimes} reports the execution times needed for computing all Bayesian bandwidth selection methods related to only one replication of  sample sizes $n=10$, 25, 50, 100, 200 and 500, and 1000 for the target function A. The results  indicate  similar performances of
Bayesian adaptive approach with  standard and modified gamma kernel smoothers. The difference in Central Processing Unit (CPU) times goes slightly to the advantage of the standard gamma  as the sample size increases.
\begin{table}[!htbp]
	\begin{center}
		\caption{Comparison of execution times (in seconds) for one replication of both  Bayesian adaptive bandwidth selections using density A.} 
		{	\begin{tabular}{rrrcccc}
				\toprule
				$n$ &\multicolumn{1}{c}{$t_{gamma}$} &\multicolumn{1}{c}{$t_{Mgamma}$}\\\midrule
				
				10 &  0.00  & 0.00  \\
				25 & 0.12  & 0.11 \\
				50 &  0.12  & 0.12\\
				100& 0.15& 0.16\\
				200 &  0.26 &  0.38 \\
				500 &  1.07 & 1.64 \\
				1000 & 3.76  & 5.96 \\
				\bottomrule
		\end{tabular}}
		\label{computtimes}
	\end{center}
\end{table}

\begin{figure}[!htbp]
		\centering
		\vspace{-1.4cm}
		\mbox{
			\stackunder{	\resizebox*{6.cm}{!}{\includegraphics{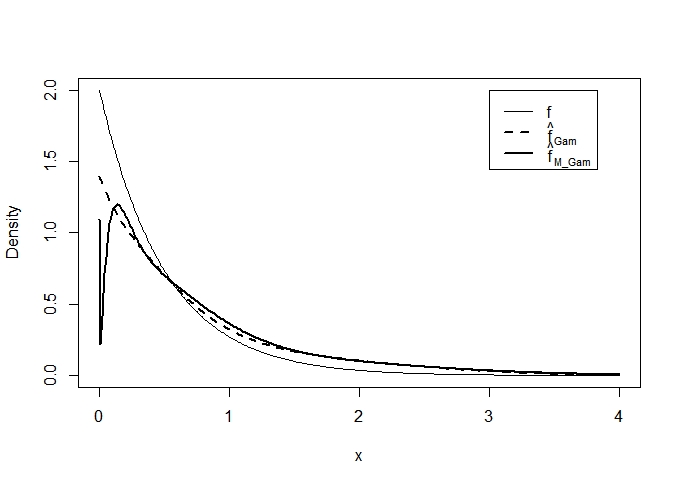}} }{(A1)}\hspace{3pt}
			\stackunder{\resizebox*{6.cm}{!}{\includegraphics{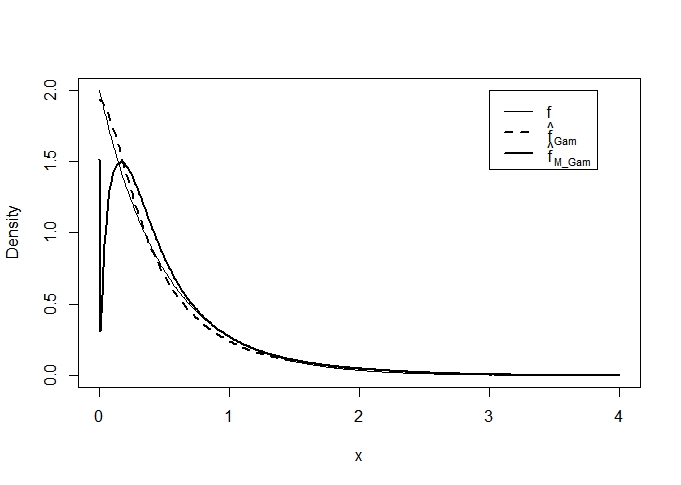}} }{(A2)}
		}\vspace{-0.1cm}
		\mbox{
				\stackunder{	\resizebox*{6.cm}{!}{\includegraphics{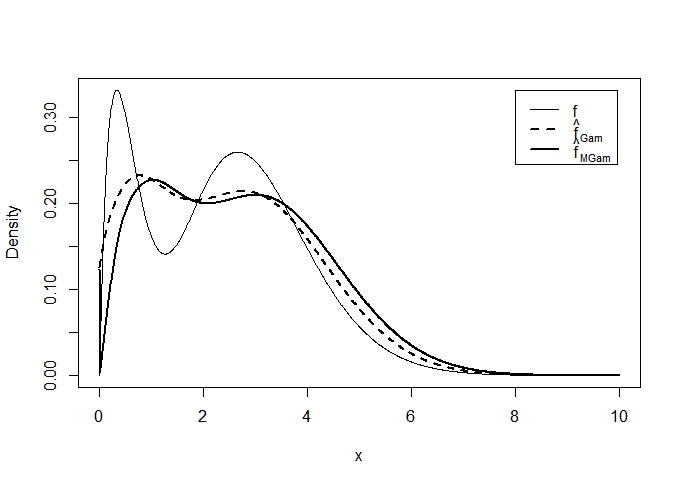}} }{(B1)}\hspace{3pt}
			\stackunder{\resizebox*{6.cm}{!}{\includegraphics{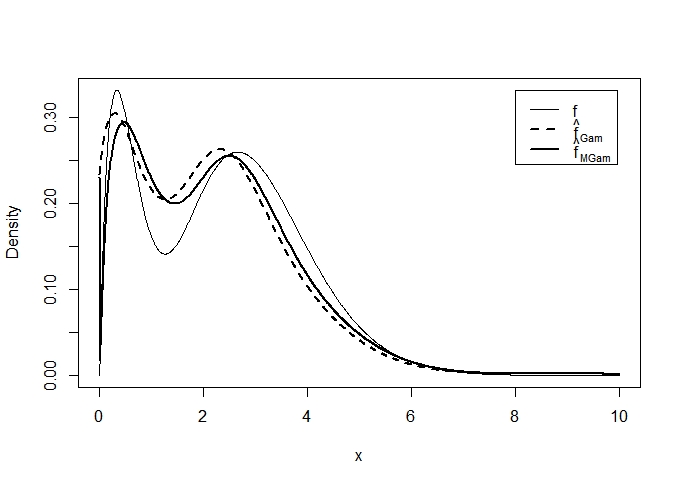}} }{(B2)} 
		}\vspace{-0.1cm}
		\mbox{
		\stackunder{	\resizebox*{6.cm}{!}{\includegraphics{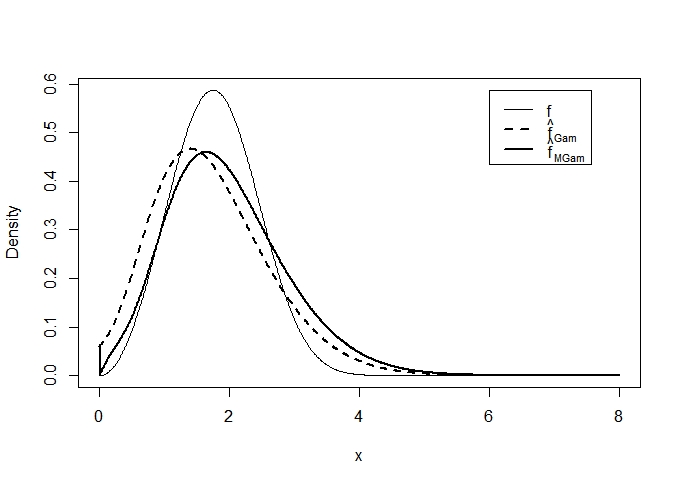}} }{(C1)}\hspace{3pt}
		\stackunder{\resizebox*{6.cm}{!}{\includegraphics{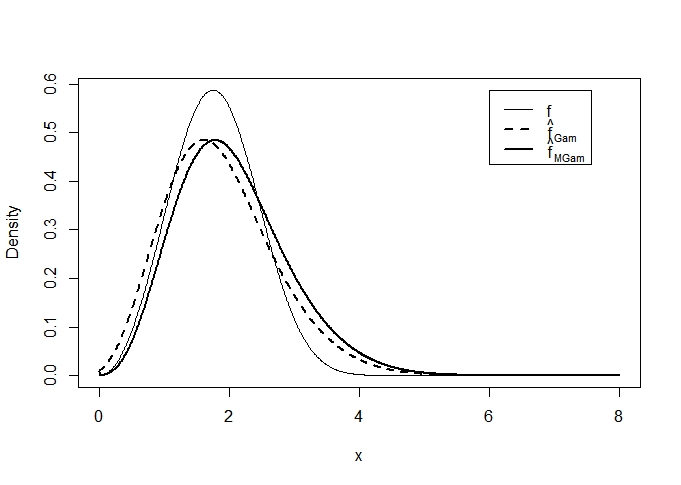}} }{(C2)} 
	}\vspace{-0.1cm}
	\mbox{
	\stackunder{	\resizebox*{6.cm}{!}{\includegraphics{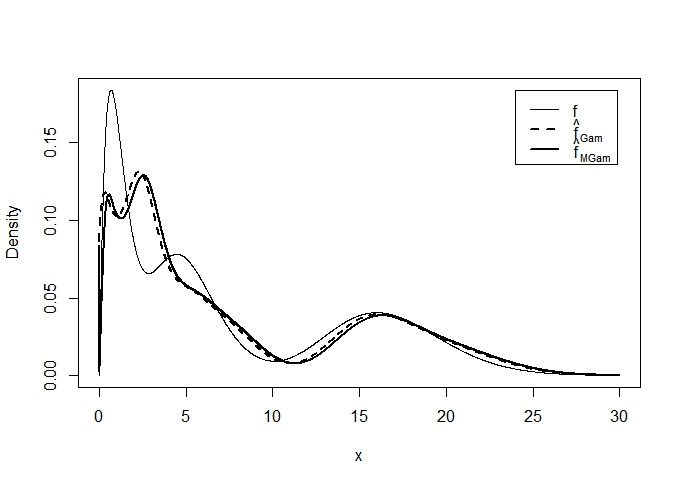}} }{(D1)}\hspace{3pt}
	\stackunder{\resizebox*{6.cm}{!}{\includegraphics{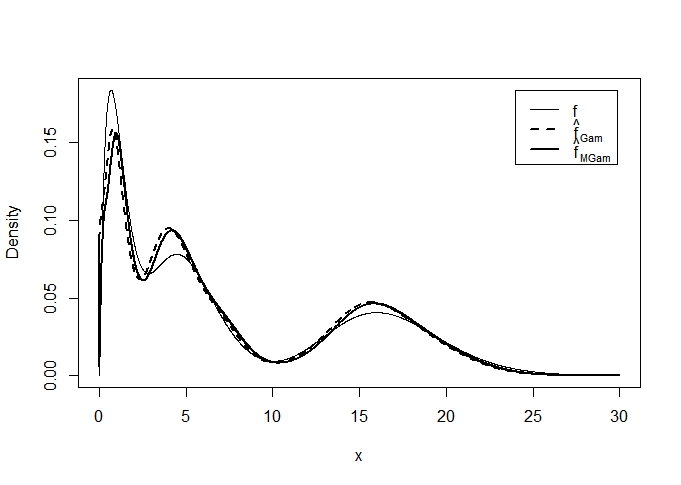}} }{(D2)}
}
	\caption{True pdf and gamma kernel smoothers on full support for A, B, C and D defined in \eqref{gammaestimator} for $d=1$ with $n=50$ (left) and
		$n=100$ (right).}\label{Est_univ}
	\end{figure}

Fig.~\ref{Est_univ} gives the true density and the smoothing densities for both gamma and modified gamma estimators with Bayesian adaptive bandwidths  related  to the considered models, and for only one replication. The graphs show that both gamma estimators have similar performances with more difficulties, for the modified one, to give satisfying estimates at the small boundary region (e.g. Part A1--A2 of Fig. \ref{Est_univ}).


Table~\ref{ISE_univ} shows some expected values of ISE with respect to the four scenarios A, B, C and D, and according to the sample sizes $n = 10, 25, 50, 100, 200, 500, 1000$.  Thus, we observe the following behaviors. The smoothings are better as the sample sizes increase for both methods. Globally, the standard gamma kernel is preferable to the modified one especially for convex densities with maximum probability appearance at origin (e.g. Part A of Table \ref{ISE_univ}). Conversely, the modified gamma suits more for unimodal densities and for all sample sizes.

Finally, Table \ref{ISE_univ} and Fig. \ref{Est_univ} indicate that the modified gamma kernel is preferable for unimodal densities with mode far from the boundary region while  the standard gamma take the lead in other situations. Then, the choice of the appropriate multiple kernel between the standard, modified and combined gamma ones arises.
\begin{table}[!ht]
	
	\begin{center}
		\caption{ Expected values ($\times 10^3$)  of $\widehat{ISE}$ and their standard deviations in parentheses  with $N=100$ replications  using	both Bayesian adaptive bandwidths for scenarios A, B, C and D. } \label{ISE_univ}
		\begin{tabular}{rrrrrrrr}
			\toprule
			&$n$&\multicolumn{1}{c}{$\widehat{ISE}_{Gamma}$ }&\multicolumn{1}{c}{$\widehat{ISE}_{Mgamma}$ }    \\\hline

			\multirow{7}{*}{A}& \multirow{1}{*}{10}   &87.53 (71.10)&281.52 (134.73)\\  
			&\multirow{1}{*}{25} &38.90 (35.44)&  214.55 (287.91)\\		  
			&\multirow{1}{*}{50}  &21.83 (16.64)&   201.53 (380.13)\\
			&\multirow{1}{*}{100} &18.96 (14.31) & 66.78 (53.07) \\
			&\multirow{1}{*}{200} & 11.57 (8.96)&  46.90 (52.29) \\
			&\multirow{1}{*}{500}& 5.44 (4.60)& 18.78 (20.21) \\
			&\multirow{1}{*}{1000}& 3.22 (2.24)&9.08 (11.83) \\  

			\hline 
			\multirow{7}{*}{B}& \multirow{1}{*}{10}  &20.23 (6.57) & 35.69 (12.15)\\		
			  	&\multirow{1}{*}{25}&15.36 (7.16)&20.39 (10.42)\\
			&\multirow{1}{*}{50}&11.41 (4.90) &12.74 (6.12)\\
			&\multirow{1}{*}{200}&4.11 (1.57)&4.31 (1.66) \\
			&\multirow{1}{*}{500}& 3.33 (1.18)& 3.08 (1.21) \\
			&\multirow{1}{*}{1000}&2.17 (0.72) &1.95 (0.64)		  
			\\
			\hline 
			\multirow{7}{*}{C}& \multirow{1}{*}{10}  &86.10 (32.99) &65.39 (23.90) \\ 
			& \multirow{1}{*}{25}&52.00 (18.56) &42.53 (14.00)  \\
			&\multirow{1}{*}{50}&36.42 (11.22)&27.66 (8.77) \\
			&\multirow{1}{*}{100}&23.51 (7.43) & 21.52 (6.48)\\
			&\multirow{1}{*}{200}&15.22 (5.17) &12.04 (3.13) \\
			&\multirow{1}{*}{500}&8.08 (3.29)& 6.41 (2.47)  \\
			&\multirow{1}{*}{1000}& 5.11 (1.77)& 3.99 (1.30) \\
			\hline
		
			\multirow{7}{*}{D}& \multirow{1}{*}{10}  &11.81 (7.79) &14.86 (7.61) \\
		& \multirow{1}{*}{25}&7.06 (2.94) &8.73 (4.15)  \\
		&\multirow{1}{*}{50}&5.11 (2.59)&5.43 (2.72) \\
		&\multirow{1}{*}{100} &2.90 (1.41)& 3.14 (1.65)\\
		&\multirow{1}{*}{200}& 2.03 (0.68)& 2.18 (0.89)\\
		&\multirow{1}{*}{500}&1.16 (0.44)& 1.29 (0.52) \\
		&\multirow{1}{*}{1000}& 0.72 (0.30)& 0.77 (0.31) \\
			\bottomrule
		\end{tabular}

	\end{center}
\end{table} 

\subsection{Bivariate and multivariate case studies}\label{SSection_Biv}
\subsubsection{Bivariate  case study}

\begin{figure}[!htbp]
	\vspace*{-1cm}
	\mbox{
		
		\subfloat[(I$_C$)]{	\resizebox*{6.5cm}{!}{\includegraphics{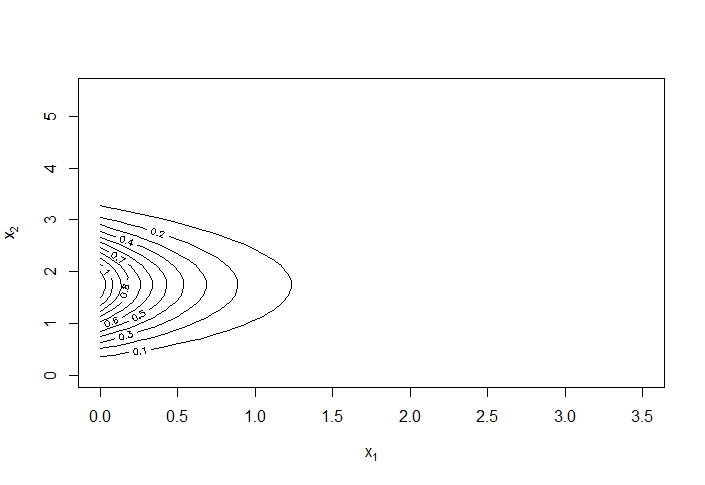}} }\hspace{3pt}
		\subfloat[(I$_S$) ]{\resizebox*{6.5cm}{!}{\includegraphics{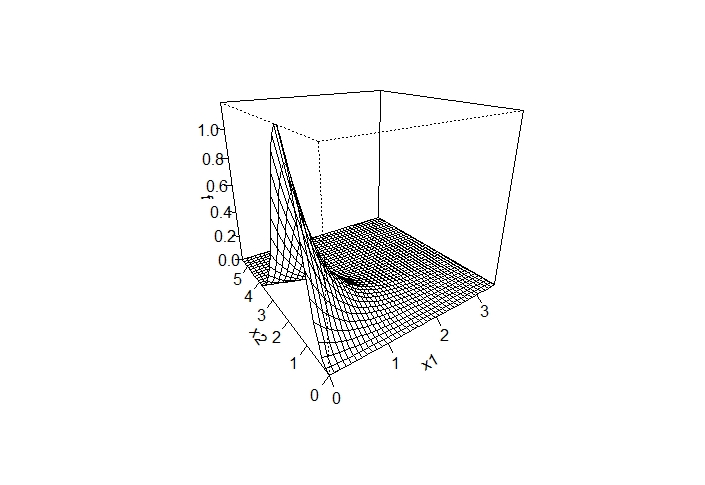}} }	
	}
	
	\mbox{
		
		\subfloat[(I$_C$-gamma)]{	\resizebox*{6.3cm}{!}{\includegraphics{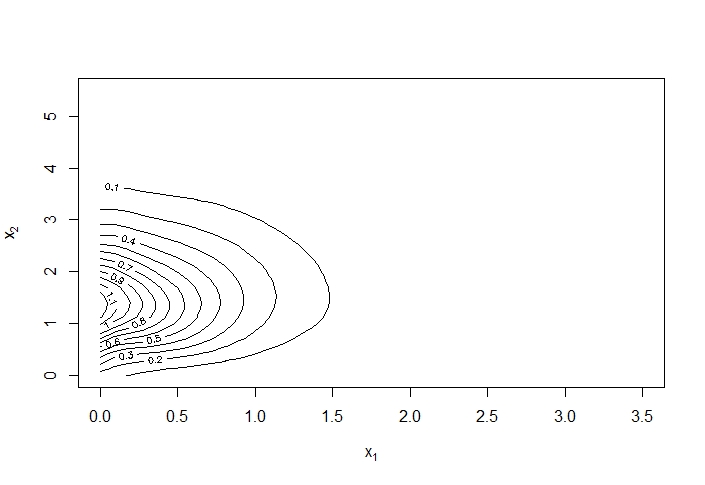}} }\hspace{3pt}
		\subfloat[(I$_S$-gamma) ]{\resizebox*{6.3cm}{!}{\includegraphics{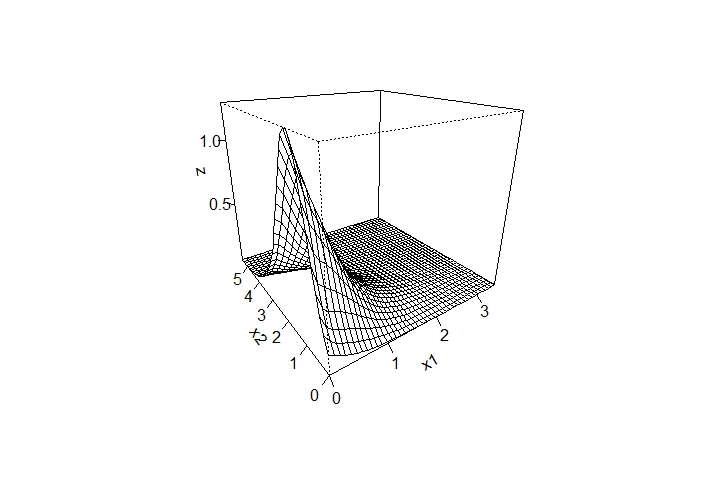}}} 
		
	}
	\mbox{ 
		
		\subfloat[(I$_C$-modified gamma)]{	\resizebox*{6.3cm}{!}{\includegraphics{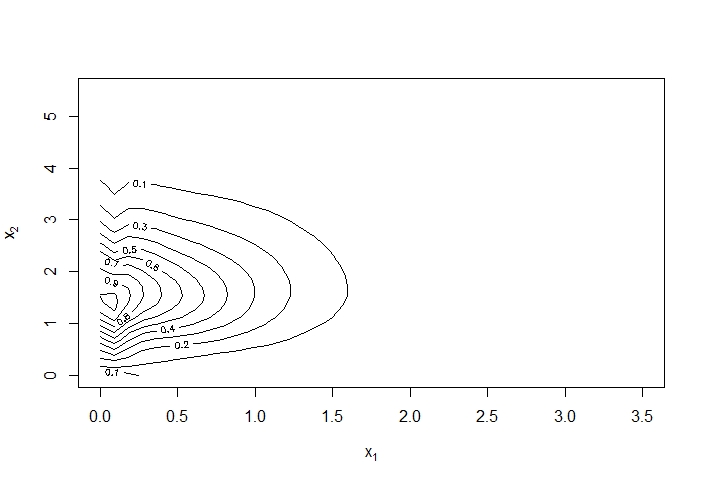}} }\hspace{3pt}
		\subfloat[(I$_S$-modified gamma)]{\resizebox*{6.3cm}{!}{\includegraphics{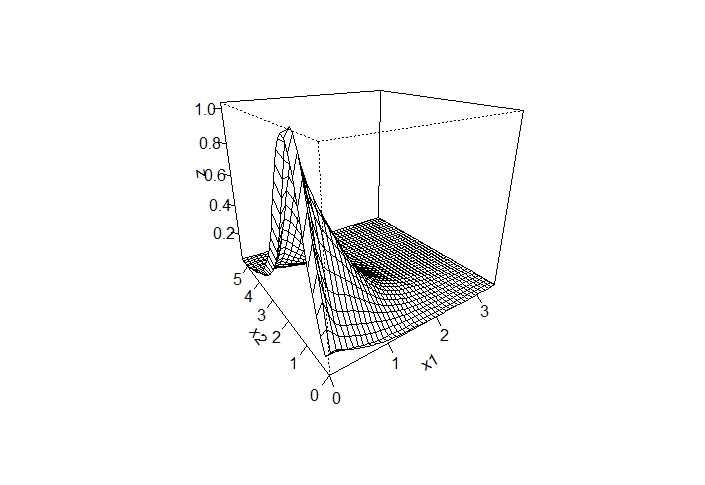}} }
	}

	\mbox{ 
	
	\subfloat[(I$_C$-combined gamma)]{	\resizebox*{6.3cm}{!}{\includegraphics{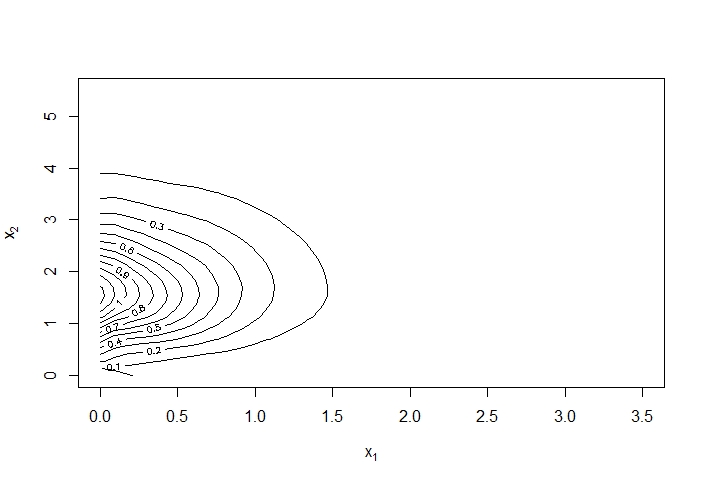}} }\hspace{3pt}
	\subfloat[(I$_S$-combined gamma)]{\resizebox*{6.3cm}{!}{\includegraphics{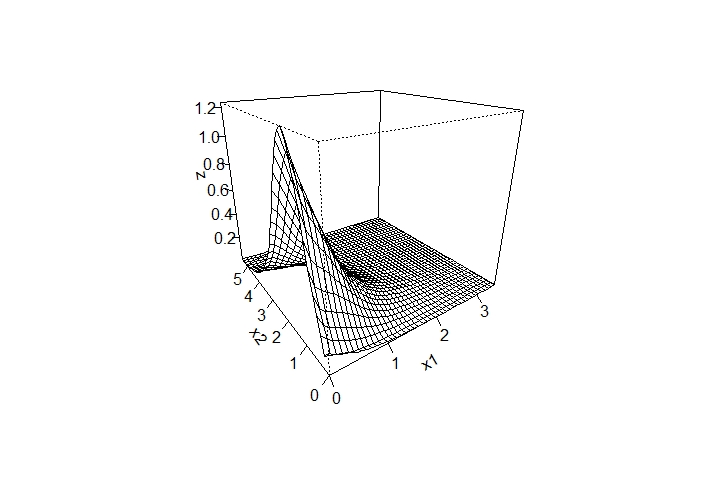}} }
}
	\caption{Contour (left) and surface (right) plots of estimated bivariate gamma-Weibull model I  according to all three Bayes selectors of bandwidths  vector $\mathbf{h}$ with $n=100$.}
	\label{dens_I_estim}
\end{figure}

\begin{figure}[!htbp]
	\mbox{
		
		\stackunder{	\resizebox*{7cm}{!}{\includegraphics{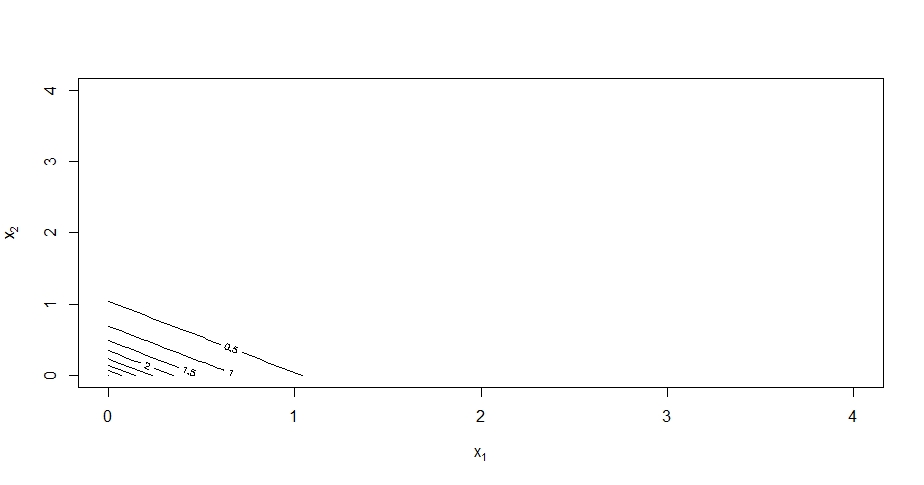}} }{(E$_C$)}\hspace{.1pt}
		\stackunder{\resizebox*{8.5cm}{!}{\includegraphics{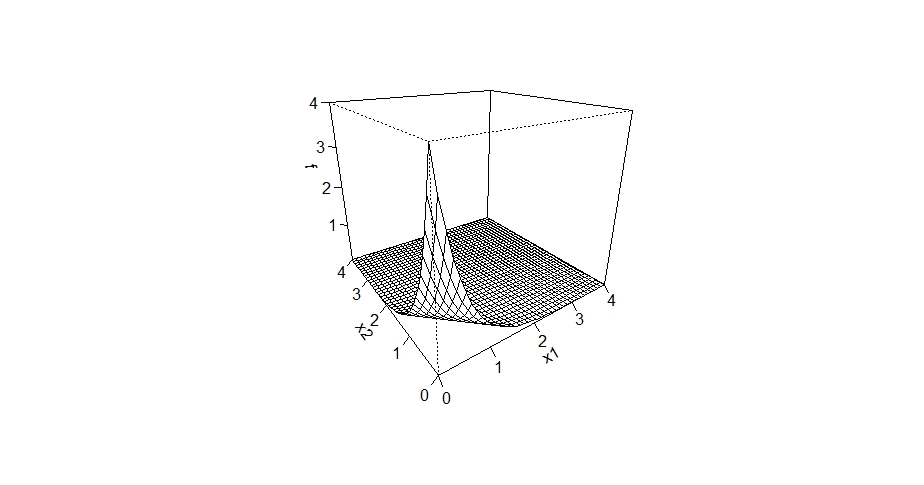}} }{(E$_S$)}
		
	}
	\mbox{
		
		\stackunder{	\resizebox*{7cm}{!}{\includegraphics{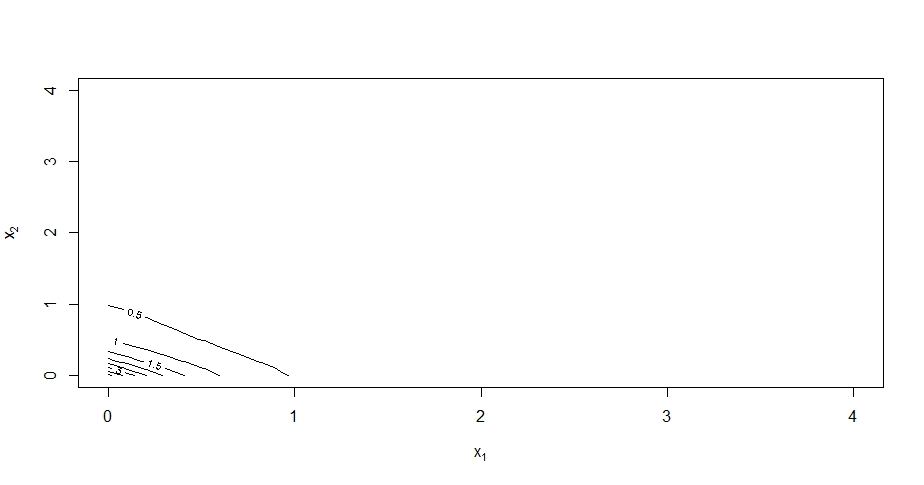}} }{(E$_C$-gamma)}\hspace{0.1pt}
		\stackunder{\resizebox*{8.5cm}{!}{\includegraphics{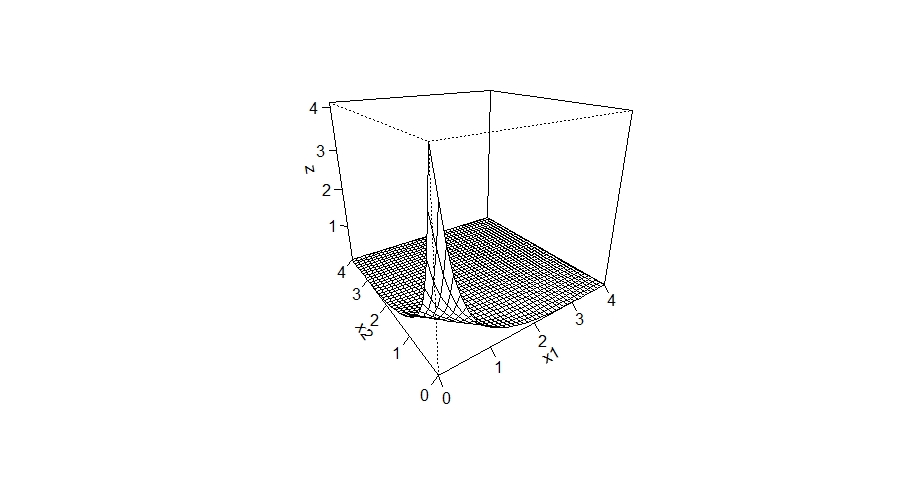}} }{(E$_S$-gamma) }
		
	}
	\mbox{ 
		
		\stackunder{	\resizebox*{7cm}{!}{\includegraphics{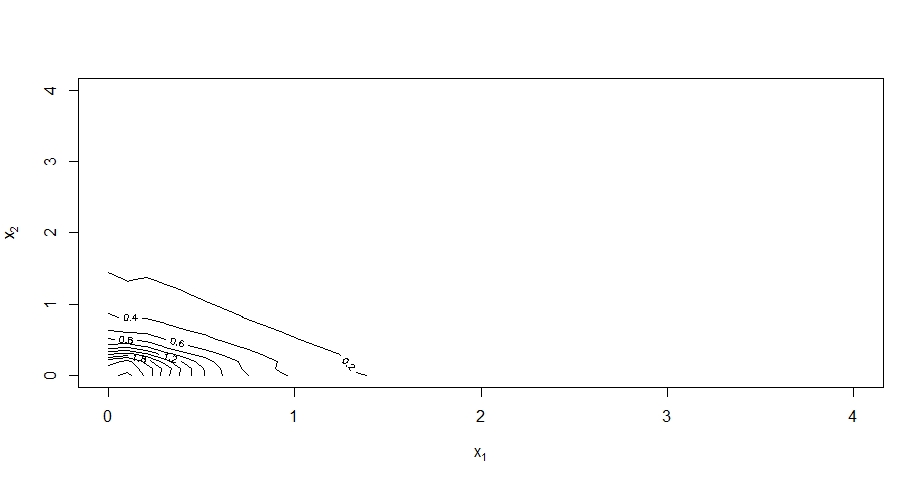}} }{(E$_C$-modified gamma)}\hspace{0.1pt}
		\stackunder{\resizebox*{8.5cm}{!}{\includegraphics{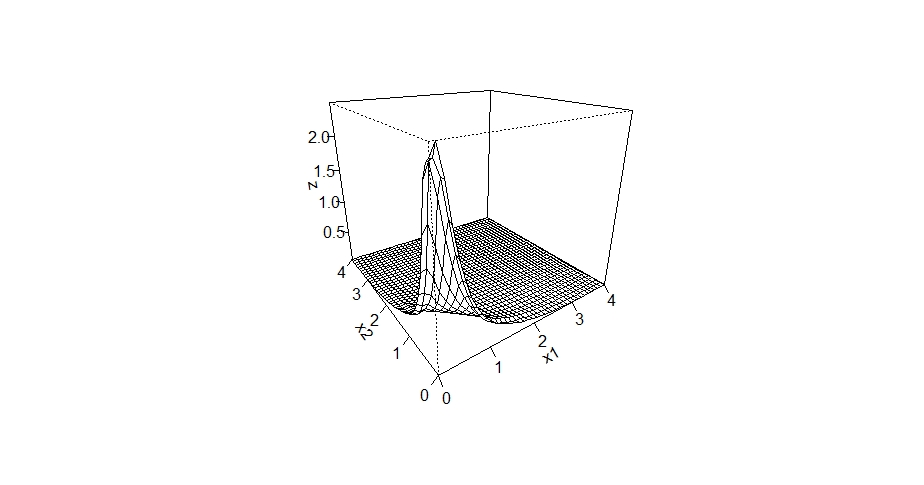}} }{(E$_S$-modified gamma)}
	}
	\mbox{ 
		
		\stackunder{	\resizebox*{7cm}{!}{\includegraphics{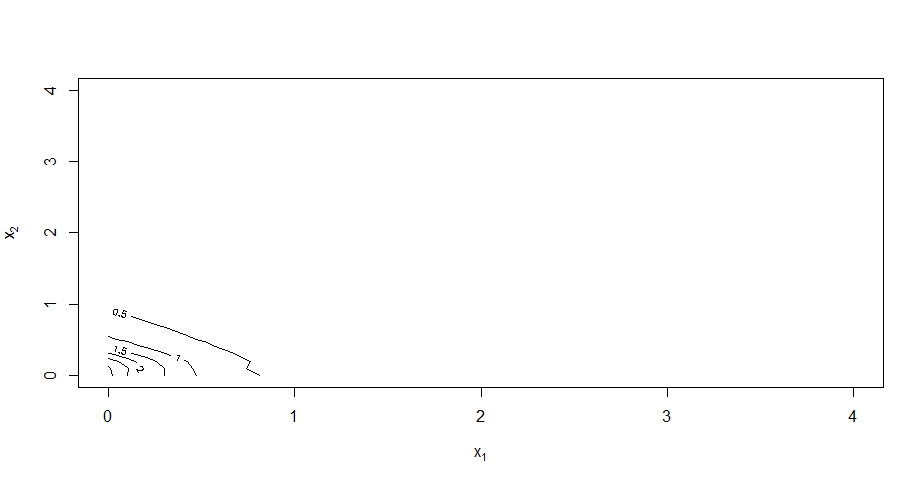}} }{(E$_C$-combined gamma)}\hspace{.1pt}
		\stackunder{\resizebox*{8.5cm}{!}{\includegraphics{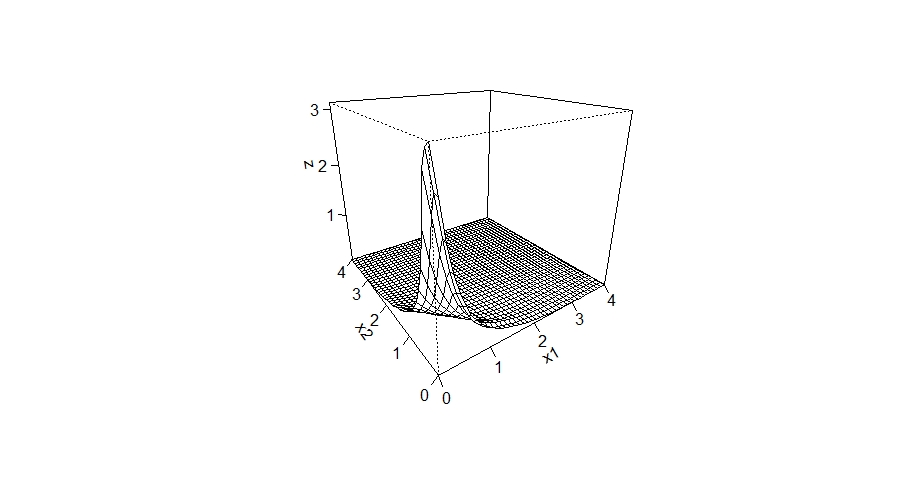}} }{(E$_S$-combined gamma)}
	}
	\caption{Contour (left) and surface (right) plots of estimated bivariate gamma model E    according to all three Bayes selectors of bandwidths vector $\mathbf{h}$ with $n=100$.}
	\label{dens_E_estim}
\end{figure}

\begin{figure}[!htbp]
		\vspace*{-2cm}
	\mbox{
		
		\stackunder{	\resizebox*{7cm}{!}{\includegraphics{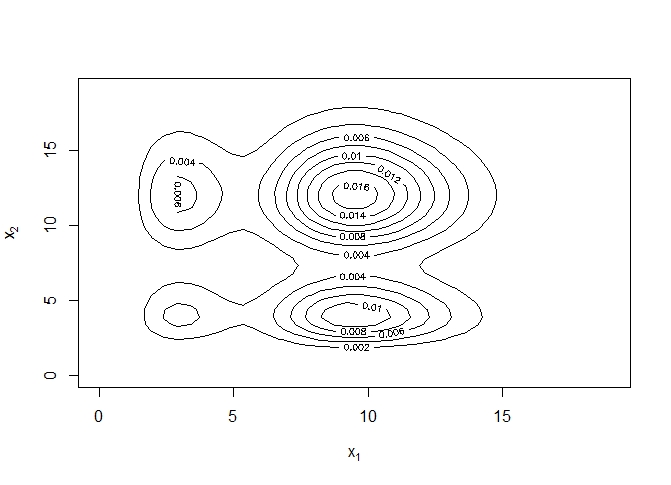}} }{(F$_C$)}\hspace{3pt}
		\stackunder{\resizebox*{7cm}{!}{\includegraphics{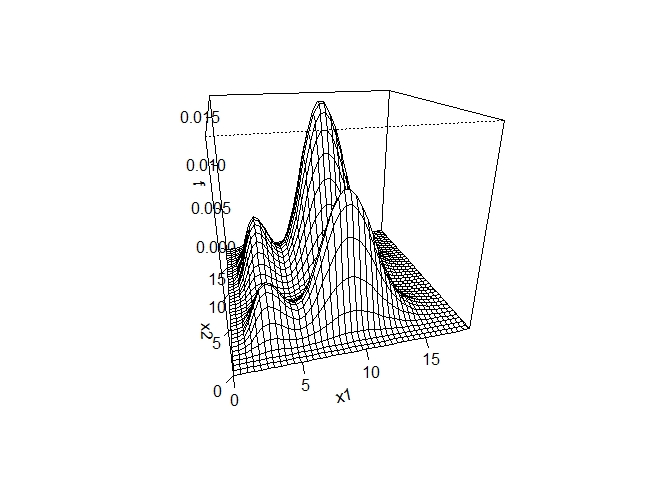}} }{(F$_S$) }
		
	}
	\mbox{
		
		\stackunder{	\resizebox*{7cm}{!}{\includegraphics{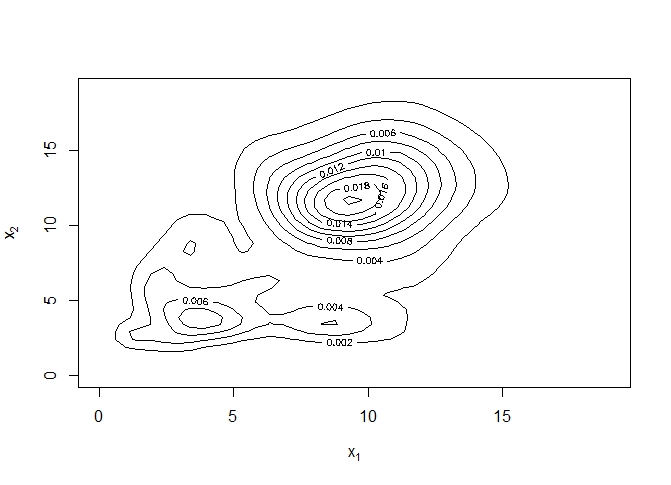}} }{(F$_C$-gamma)}\hspace{3pt}
		\stackunder{\resizebox*{7cm}{!}{\includegraphics{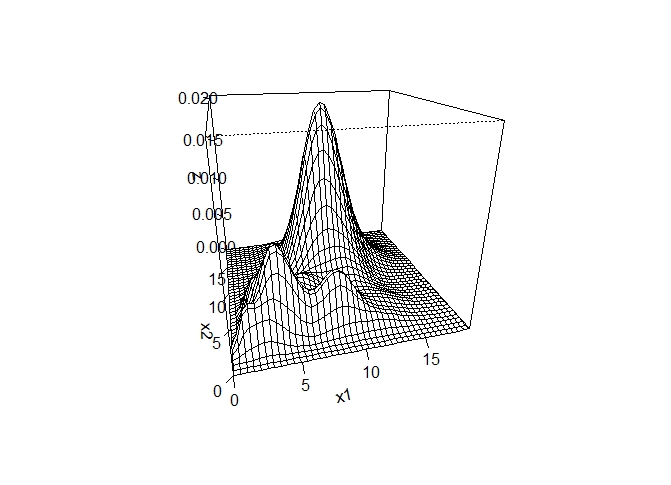}} }{(F$_S$-gamma) }
		
	}
	\mbox{ 
		
		\stackunder{	\resizebox*{7cm}{!}{\includegraphics{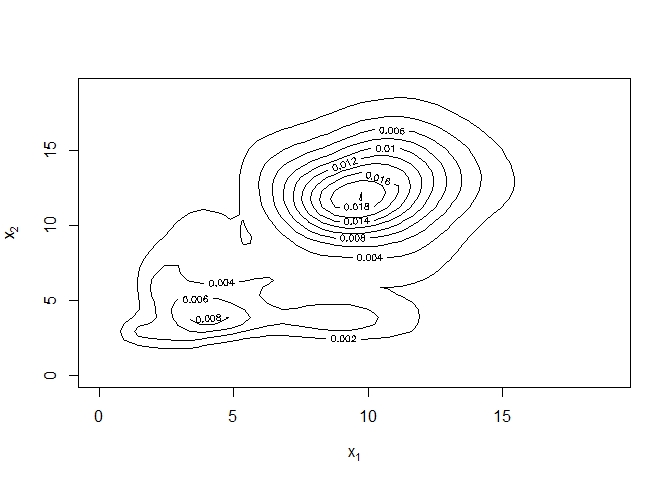}} }{(F$_C$-modified gamma)}\hspace{3pt}
		\stackunder{\resizebox*{7cm}{!}{\includegraphics{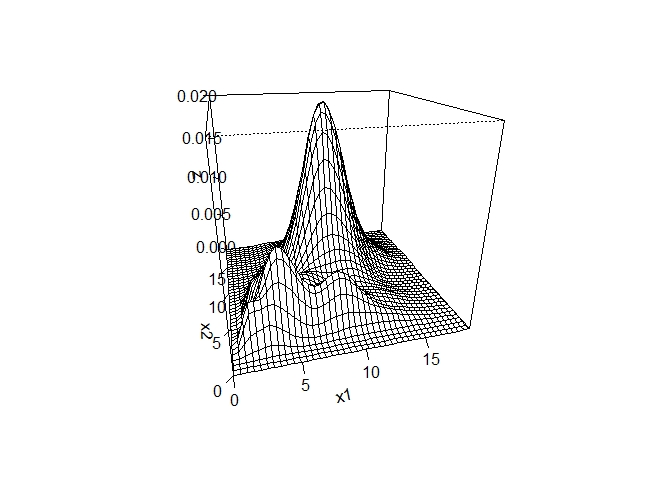}} }{(F$_S$-modified gamma)}
	}
	\mbox{ 
		
		\stackunder{	\resizebox*{7cm}{!}{\includegraphics{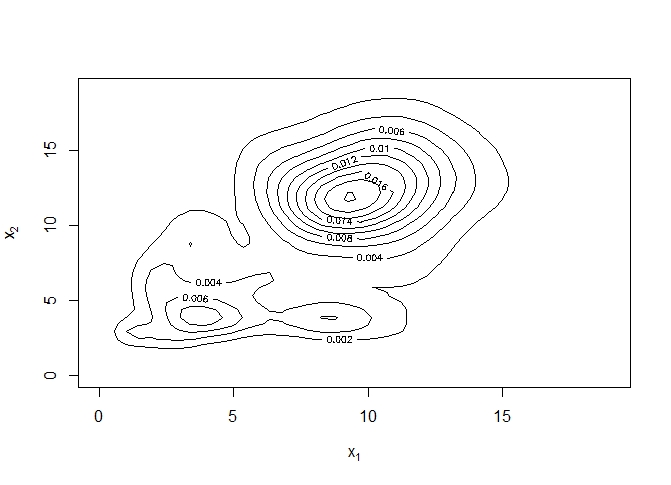}} }{(F$_C$-combined gamma)}\hspace{3pt}
		\stackunder{\resizebox*{7cm}{!}{\includegraphics{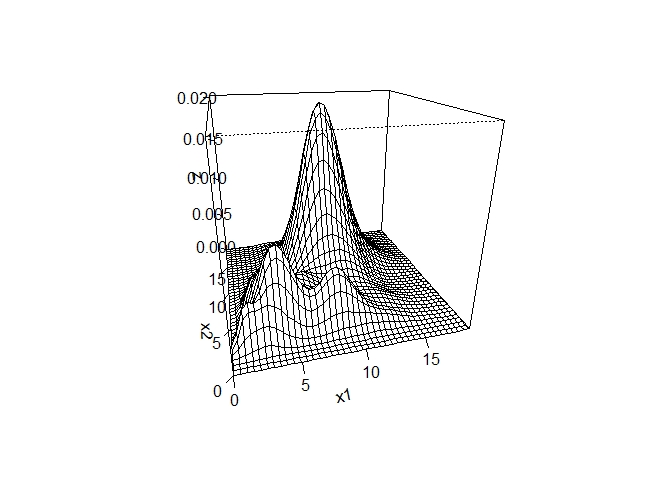}} }{(F$_S$-combined gamma)}
	}
	
	\caption{Contour (left) and surface (right) plots of estimated bivariate mixture of Erlang model F   according to all three Bayes selectors of bandwidths  vector $\mathbf{h}$ with $n=100$.}
	\label{dens_F_estim}
\end{figure}

\begin{figure}[!htbp]
	\mbox{
		
		\stackunder{	\resizebox*{7cm}{!}{\includegraphics{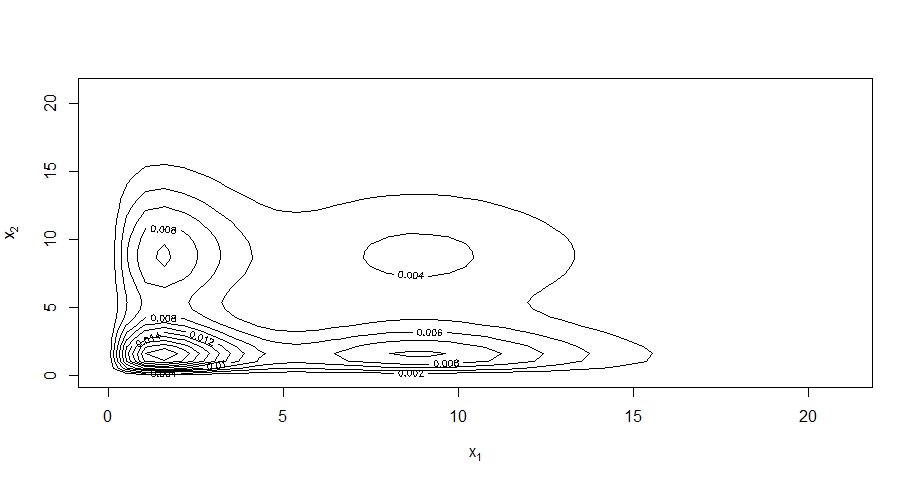}} }{(G$_C$)}\hspace{3pt}
		\stackunder{\resizebox*{8.5cm}{!}{\includegraphics{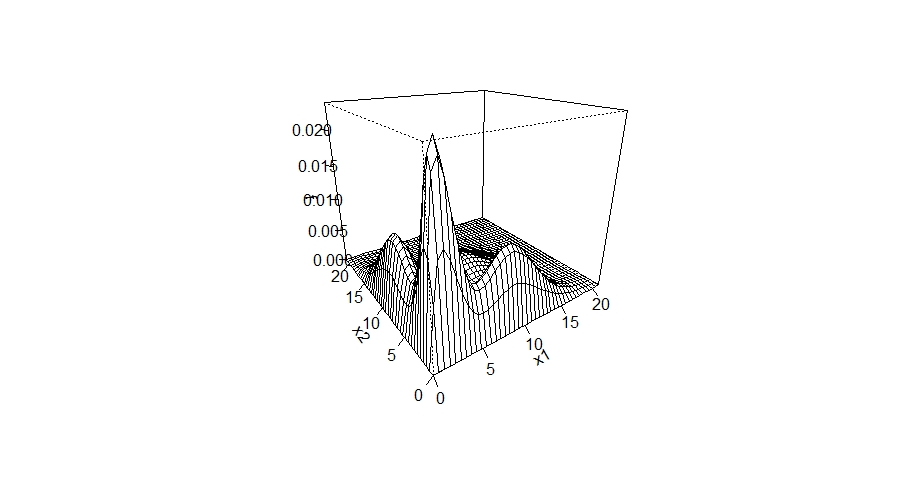}} }{(G$_S$) }
		
	}
	
	\mbox{
		
		\stackunder{	\resizebox*{7cm}{!}{\includegraphics{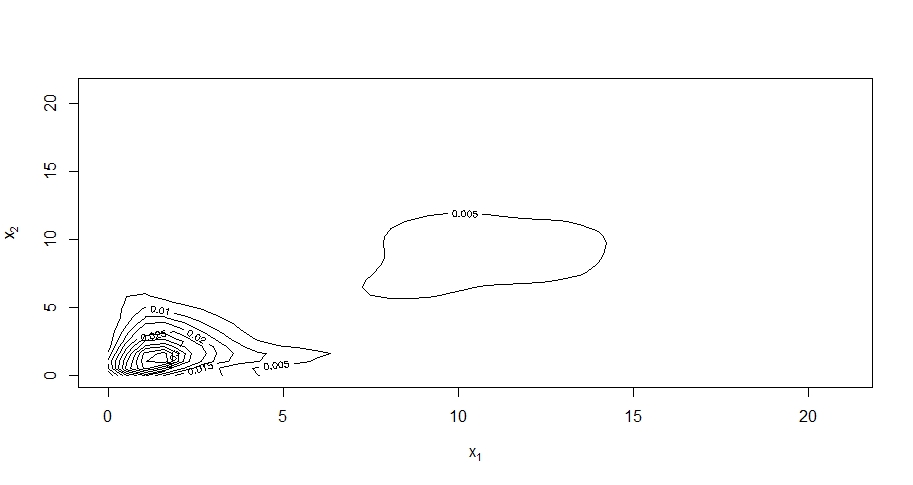}} }{(G$_C$-gamma)}\hspace{3pt}
		\stackunder{\resizebox*{8.1cm}{!}{\includegraphics{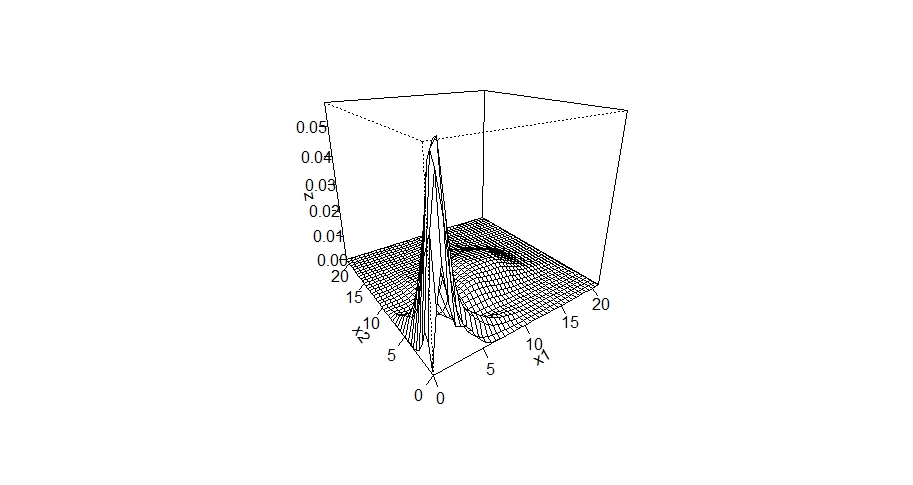}} }{(G$_S$-gamma) }
		
	}
	\mbox{ 
		
		\stackunder{	\resizebox*{7cm}{!}{\includegraphics{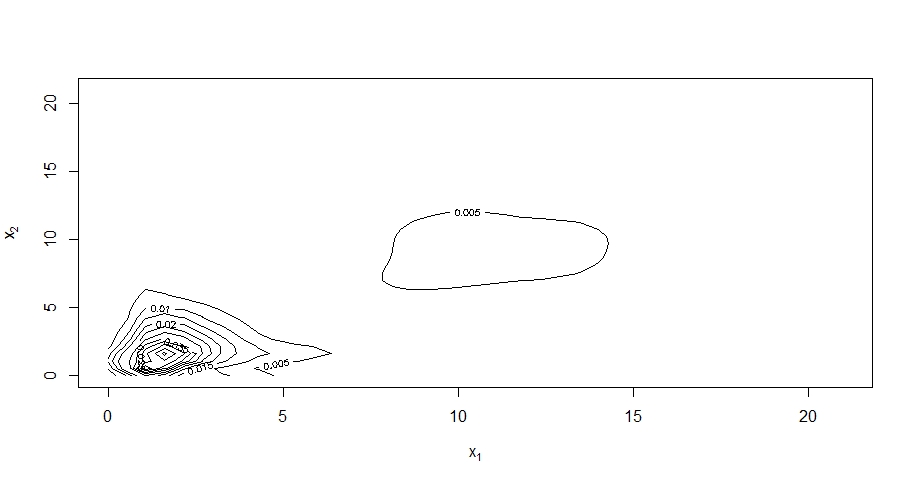}} }{(G$_C$-modified gamma)}\hspace{3pt}
		\stackunder{\resizebox*{8.1cm}{!}{\includegraphics{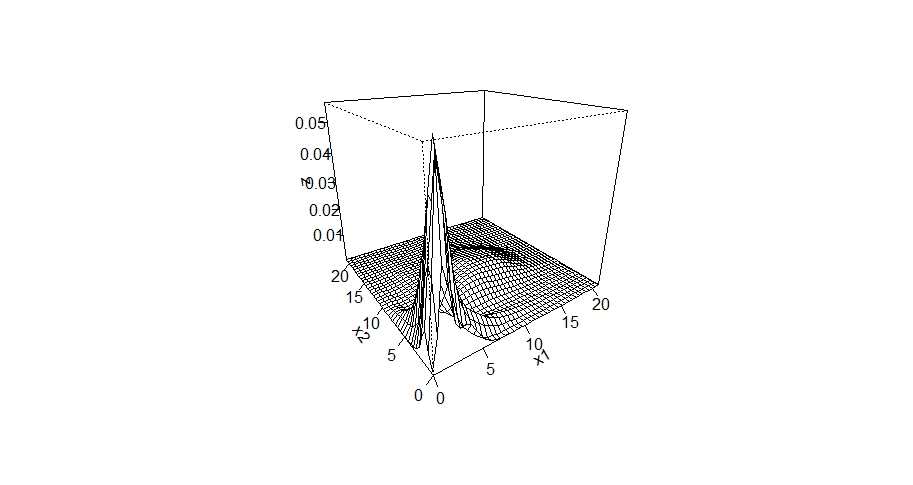}} }{(G$_S$-modified gamma)}
	}
	\mbox{ 
		
		\stackunder{	\resizebox*{7cm}{!}{\includegraphics{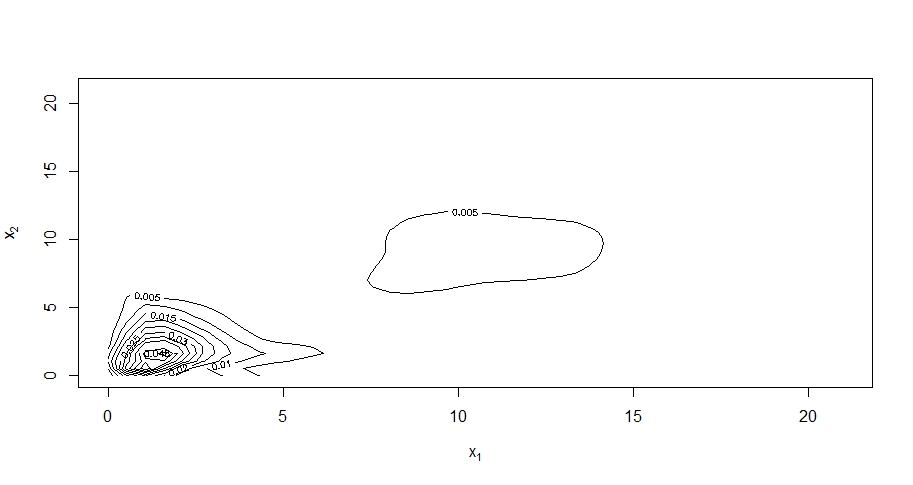}} }{(G$_C$-combined gamma)}\hspace{3pt}
		\stackunder{\resizebox*{8.1cm}{!}{\includegraphics{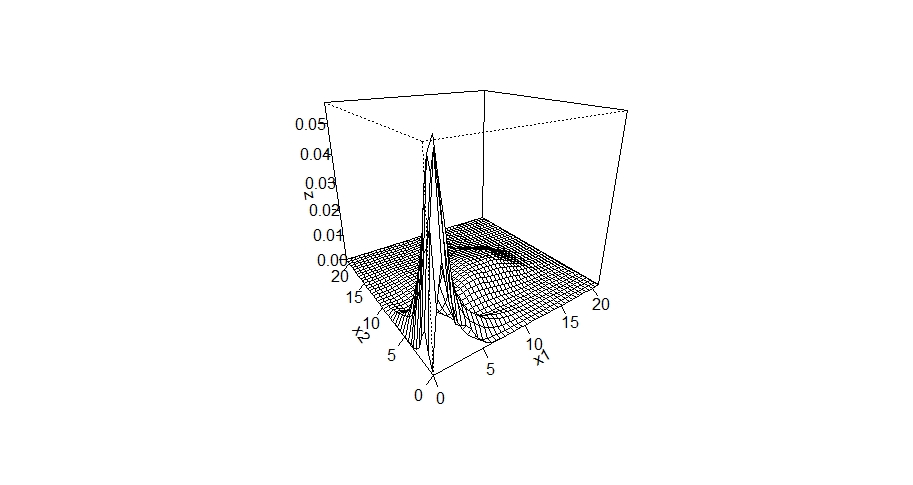}} }{(G$_S$-combined gamma)}
	}
	\caption{Contour (left) and surface (right) plots of estimated bivariate mixture of gamma model G    according to all Bayes selectors of bandwidths  vector $\mathbf{h}$ with $n=100$.}
	\label{dens_G_estim}
\end{figure}

\begin{figure}
		\vspace*{-1cm}
	\mbox{
		
		\stackunder{	\resizebox*{7cm}{!}{\includegraphics{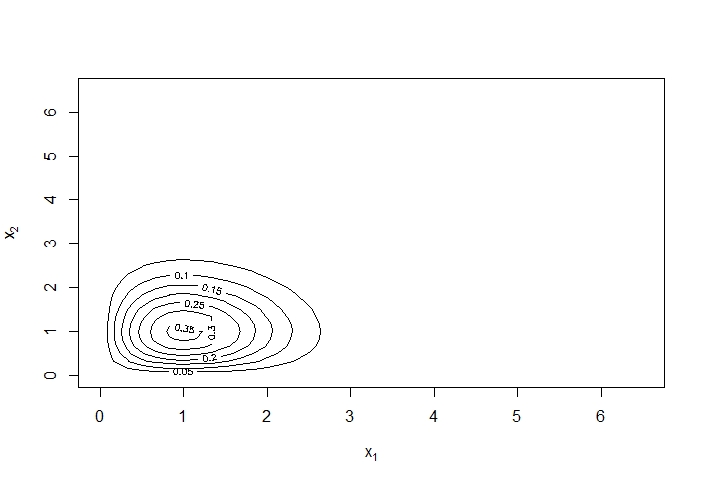}}}{(H$_C$)} \hspace{3pt}
		\stackunder{\resizebox*{7cm}{!}{\includegraphics{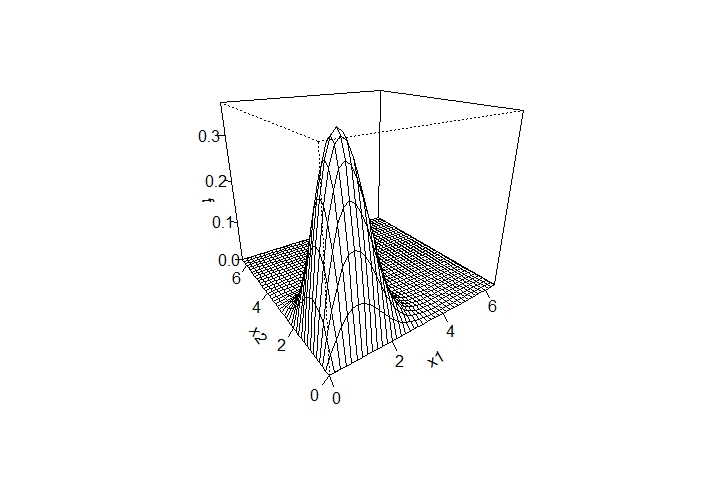}} }{(H$_S$) }	
	}
	
	\mbox{
		
		\stackunder{	\resizebox*{7cm}{!}{\includegraphics{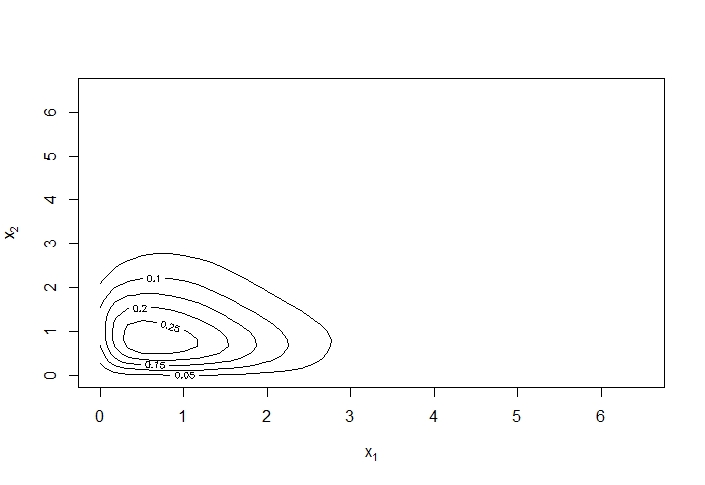}} }{(H$_C$-gamma)}\hspace{3pt}
		\stackunder{\resizebox*{7cm}{!}{\includegraphics{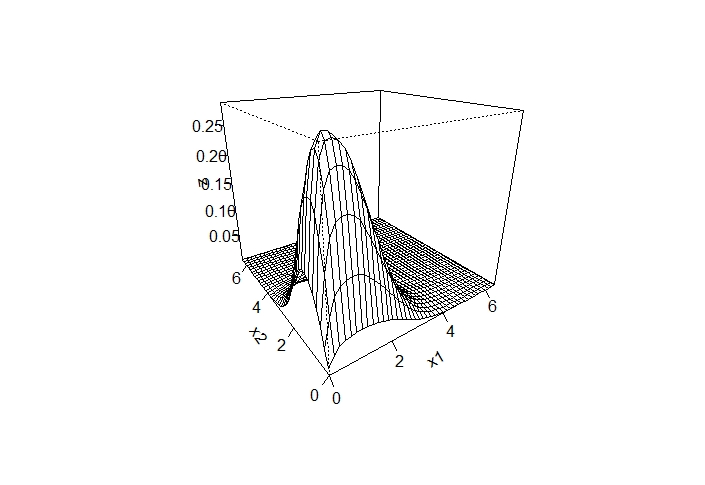}} }{(H$_S$-gamma) }
		
	}
	\mbox{ 
		
		\stackunder{	\resizebox*{7cm}{!}{\includegraphics{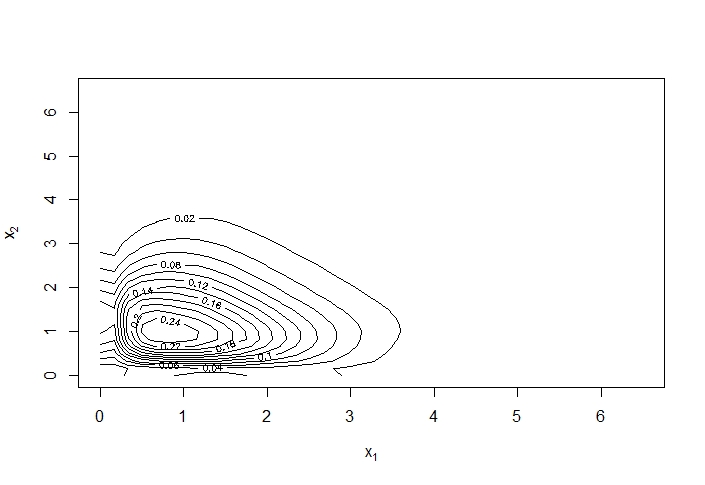}} }{(H$_C$-modified gamma)}\hspace{3pt}
		\stackunder{\resizebox*{7cm}{!}{\includegraphics{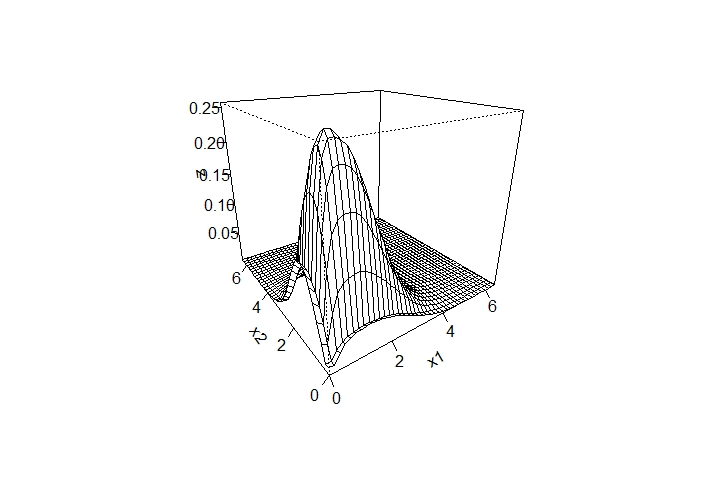}} }{(H$_S$-modified gamma)}
	}
	\mbox{ 
		
		\stackunder{	\resizebox*{7cm}{!}{\includegraphics{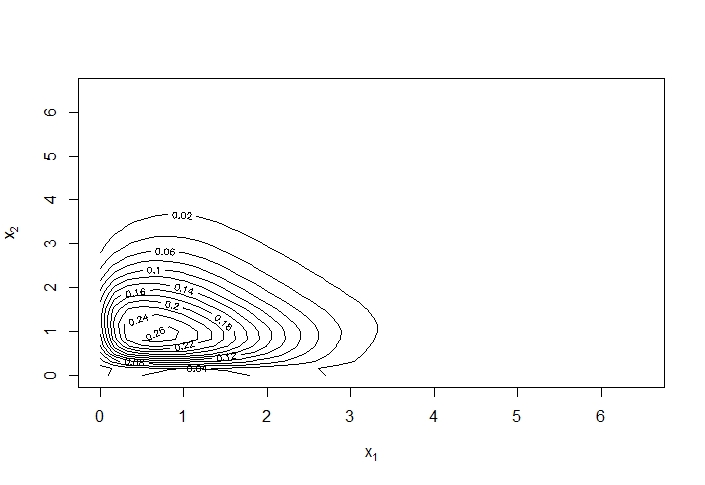}} }{(H$_C$-combined gamma)}\hspace{3pt}
		\stackunder{\resizebox*{7cm}{!}{\includegraphics{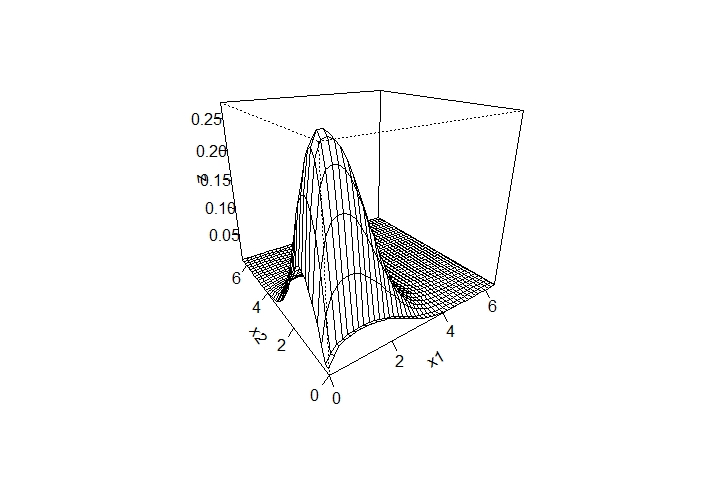}} }{(H$_S$-combined gamma)}
	}
	\caption{Contour (left) and surface (right) plots of estimated bivariate Weibull model H    according to all Bayes selectors of bandwidths  vector $\mathbf{h}$ with $n=100$.}
	\label{dens_H_estim}
\end{figure}

We here designed five bivariate scenarios, denoted by E, F, G, H, and I,  to simulate bivariate nonnegative datasets. These scenarios aim to justify the use of the compromise combined estimator \eqref{mix_gammaestimator}, in addition to the same motivations of the univariate case.  Both surface and contour plots for the corresponding true and smoothed densities are given in Fig.~\ref{dens_I_estim} for Scenario I,  and in Fig.~\ref{dens_E_estim}--\ref{dens_H_estim} of additional material for Scenarios E--H.
\begin{itemize}	
	\item Scenario E is generated by using the bivariate gamma density 
	$$f_{E}(x_1,x_2)=2  e^{-2 x_1}\times2  e^{-2 x_2}1_{[0,\infty)^2}(x_1,x_2);$$	
	\item Scenario F comes from a bivariate mixture of Erlang distributions 
	$$	f_F(x_1,x_2)=\left(\frac{2}{11}\frac{2^{7} x_1^{6} e^{-2 x_1}}{6 !}+\frac{9}{11}\frac{2^{20} x_1^{19} e^{-2 x_1}}{19 !}\right)\times \left(\frac{2}{7}\frac{2^{9} x_2^{8} e^{-2 x_2}}{8 !}+\frac{5}{7}\frac{2^{25} x_2^{24} e^{-2 x_2}}{24 !}\right)1_{[0,\infty)^2}(x_1,x_2);$$ 

	\item Scenario G is generated by using a  bivariate mixture gamma densities 
		\begin{align*}
			f_G(x_1,x_2)&=\left(0.5 \frac{x_1^{1.5} \exp (- x_1)}{\Gamma\left(2.5\right)}+0.5 \frac{x_1^{9} \exp (- x_1)}{\Gamma\left(10\right)}\right)\\\nonumber
		&\quad \times\left(0.5 \frac{x_2^{1.5} \exp (- x_2)}{\Gamma\left(2.5\right)}+0.5 \frac{x_2^{9} \exp (- x_2)}{\Gamma\left(10\right)}\right)1_{[0,\infty)^2}(x_1,x_2);
		\end{align*}		
\item Scenario H is a bivariate Rayleigh density
$$f_H(x_1,x_2)=x_1\exp(-x_1^2/2)\times  x_2\exp(-x_2^2/2)1_{[0,\infty)^2}(x_1,x_2).$$
\item Scenario I is a product of two univariate gamma and Weibull densities
$$f_I(x_1,x_2)= \left(2  e^{-2 x_1}\right)\times \frac{3}{2}\left(\frac{x_2}{2}\right)^{2} \exp\left(-[x_2/ 2]^{3}\right)1_{[0,\infty)^2}(x_1,x_2).$$
\end{itemize}

Table~\ref{computtimes1} reports the execution times needed for computing all Bayesian bandwidths selection methods related to only one replication of  sample sizes $n=10$, 25, 50, 100, 200 and 500, and for the target function E. Once again, we observe similar performance in CPU times from both multiple standard and modified gamma kernels with a greater advantage for the first  as the sample size increases. The pure combined gamma kernel is a true compromise between  the standard and modified ones for all sample sizes, and in terms of CPU times.

 \begin{table}[!ht]
	\begin{center} 	
		\caption{Comparison of execution times (in seconds) for one
			replication of all Bayesian adaptive bandwidth selections by using density E.} 
		{	\begin{tabular}{rrrcccc}
				\toprule
				$n$ &\multicolumn{1}{c}{$t_{gamma}$} &\multicolumn{1}{c}{$t_{Mgamma}$} &\multicolumn{1}{c}{$t_{Gam-Mgamma}$}\\\midrule
				20 &  0.03   & 0.02&0.00 \\
				50 & 0.13  & 0.12 	&0.12\\
				100 & 0.46 & 0.62&0.61 \\
				200& 1.56&2.45&2.37\\
				500 &10.92  & 15.92 & 14.21 \\
				1000 &  41.29&61.98&56.01 \\
				\bottomrule
		\end{tabular}}
		\label{computtimes1}
	\end{center}
\end{table}

Fig.~\ref{dens_I_estim} displays the surface and contour plots of the true and smoothed densities with  Bayesian adaptive bandwidths for the  model I, and for one replication; see also the graphs in supplemental material for  models E, F, G and H. Once again, the graphs shed light  the good smoothing performance of all methods in general.

Table~\ref{ISE_biv} shows some expected values of $\widehat{ISE}$  with respect  to the five scenarios E, F, G, H and I, and according  to the sample sizes $n=20, 50, 100, 200, 500$. Then, we  observe the following behaviors. The smoothings are improving for all methods as the sample size increases. In terms of performance, the combined gamma kernel is a good compromise between both standard and gamma kernels (e.g. Parts E--F of Table~\ref{ISE_biv}). It is even the best of the three methods when the multiple kernel is build from the most performing univariate  (standard or modified) kernels according to each univariate margin (e.g. I of Table~\ref{ISE_biv}). As in the univariate cases, the pure modified gamma kernel is preferable for unimodal densities (e.g. Part G of Table~\ref{ISE_biv}) whereas the multiple standard kernel suit more for convex densities and also for multimodal densities with modes in the interior region of the support (e.g. part E and F of  Table~\ref{ISE_biv}). For Part H of  Table~\ref{ISE_biv}, the combined gamma kernel is the best for small sample sizes while  the pure multiple gamma kernel is suitable for moderate and larger ones.

\begin{table}[!ht]
		
	\begin{center}
		\caption{ Expected values ($\times 10^3$) and their standard deviations in parentheses of $\widehat{ISE}$ with $N=100$ replications  using Bayesian adaptive bandwidths for all bivariate scenarios E, F, G, H and I.} 
	{	\begin{tabular}{rrrrrrr}
		\toprule
			&$n$&\multicolumn{1}{c}{$\widehat{ISE}_{Gamma}$}&\multicolumn{1}{c}{$\widehat{ISE}_{Mgamma}$ } &\multicolumn{1}{c}{$\widehat{ISE}_{Gam-Mgam}$ }  \\\midrule		
			
			\multirow{5}{*}{E}& \multirow{1}{*}{20} &99.72 (55.21)&328.26 (156.21) & 240.36 (124.71)\\
			&\multirow{1}{*}{50} &48.68 (33.78)  &153.06 (74.64)& 103.24 (57.24) \\
			&\multirow{1}{*}{100}&40.85 (24.95) &97.11 (48.24) & 68.29 (33.35)
			  \\  
				&\multirow{1}{*}{200}&25.25 (16.86) & 55.84 (30.98) &39.85 (23.05)
				  \\				    
			&\multirow{1}{*}{500}&14.69 (6.90)& 28.75 (12.20) & 21.34 (9.38)\\

			\midrule 
			\multirow{5}{*}{F}& \multirow{1}{*}{20}  &2.50 (0.90)&2.55 (0.95) &2.51 (0.83)\\
				&\multirow{1}{*}{50}&2.09 (0.54)&2.07 (0.60)&2.05 (0.56)
			\\
			&\multirow{1}{*}{100}&1.93 (0.56)&1.92 (0.54)&1.92 (0.55)\\
			&\multirow{1}{*}{200}&1.82 (0.56)&1.79 (0.49)&1.81 (0.54)\\
			&\multirow{1}{*}{500}&1.71 (0.14)&1.70 (0.14) & 1.71 (0.14)
			   \\
			\midrule 

			\multirow{5}{*}{G}& \multirow{1}{*}{20}  &6.04 (2.18)& 5.79 (2.42) & 5.77 (2.00) 
			   \\		
			&\multirow{1}{*}{50}&5.38 (1.54)&4.76 (1.29) &5.15 (1.46)  \\
			&\multirow{1}{*}{100}&5.46 (1.07)& 4.85 (0.90) & 5.14 (1.00)		  
			   \\
			&\multirow{1}{*}{200}&5.22 (0.62)&4.86 (0.55)&5.01 (0.57)\\
			&\multirow{1}{*}{500}&5.14 (0.47)& 4.88 (0.45) &5.01 (0.46)
			  \\
				\midrule
			
			\multirow{5}{*}{H}& \multirow{1}{*}{20}  & 45.04 (16.66)& 47.57 (16.21) & 44.94 (15.13)
			   \\			    		 
			&\multirow{1}{*}{50}&28.12 (8.27)& 28.08 (8.02) & 27.54 (8.48)			  
 \\
			&\multirow{1}{*}{100}&20.19 (6.18)& 19.46 (5.68) & 19.53 (6.10)
	 \\
			&\multirow{1}{*}{200}&11.85 (4.12)&11.60 (2.89) & 12.10 (3.12)  
		   \\
			&\multirow{1}{*}{500}&7.06 (2.15)&6.65 (2.03)&6.77 (1.88)
			   \\	
			   	  \midrule
			   	  
			   \multirow{5}{*}{I}& \multirow{1}{*}{20}  &74.36 (27.60)&99.50 (106.35)&62.00 (18.19)\\
			   &\multirow{1}{*}{50}&53.19 (18.31)& 71.65 (10.22) &46.44 (11.32)\\
			   &\multirow{1}{*}{100} 
			   &51.98 (14.19)& 44.53 (10.40) &43.20 (9.38)
			   \\
			  &\multirow{1}{*}{200}&51.93 (12.38)  & 43.54 (8.77) & 45.61 (9.97)\\
			   &\multirow{1}{*}{500}& 48.26 (10.30) &41.71 (7.44) & 42.83 (9.04)
			   \\	
			\bottomrule
		\end{tabular}}
		\label{ISE_biv}
	\end{center}
	
\end{table} 

\subsubsection{Multivariate case study}

In order to investigate simulations for higher dimensions $d>2$, three test pdfs of Scenarios K, L and M are considered for $d=3$ and 5, respectively. 
\begin{itemize}
	\item Scenario K is generated by a 3-variate Weibull density 
		$$f_K(x_1,x_2,x_3)=\prod_{\ell=1}^3\left(\frac{x_\ell}{2}\right) \mathrm{e}^{-(x_\ell / 2)^{2}}1_{[0,\infty)^3}(x_1,x_2,x_3);$$ 
	\item Scenario L is from a 3-variate gamma density
	$$f_{L}(x_1,x_2,x_3)=\prod_{\ell=1}^3 \frac{x_\ell^{2} \exp (- x_\ell/2)}{2^3\Gamma\left(3\right)}1_{[0,\infty)^3}(x_1,x_2,x_3);$$
		\item Scenario M is a mixture of two 5-variate gamma densities
	$$f_{M}(x_1,\ldots,x_5)=\prod_{\ell=1}^5 \left(0.5 \frac{x_\ell \exp (- x_\ell/2)}{\Gamma\left(2\right)}+0.5 \frac{x_\ell^{7} \exp (- x_\ell/2)}{\Gamma\left(8\right)}\right)1_{[0,\infty)^5}(x_1,\ldots,x_5).$$
\end{itemize}

\begin{table}[!htbp]
		\caption{ Expected values ($\times 10^3$) and their standard deviations in parentheses of $\widehat{ISE}$ with $N=100$ replications  using Bayesian adaptive bandwidths for all multivariate scenarios  K, L and M.} 
		\centering
		{\begin{tabular}{rrrrrrr}
			\toprule

	&$n$&\multicolumn{1}{c}{$\widehat{ISE}_{Gamma}$}&\multicolumn{1}{c}{$\widehat{ISE}_{Mgamma}$ } &\multicolumn{1}{c}{$\widehat{ISE}_{Gam^{1}\times Mgam^{2}}$ }&\multicolumn{1}{c}{$\widehat{ISE}_{Gam^{3}\times Mgam^{2}}$ }
	   \\\midrule
	\multirow{6}{*}{K}& \multirow{1}{*}{20}   &9.32 (2.75)&30.00  (4.41) & 9.65 (2.35)&--\\ 
	&\multirow{1}{*}{50} &5.88 (2.17) &27.26   (5.07)&6.42 (2.20)&--\\
	&\multirow{1}{*}{100}& 3.89 (1.23)  &22.98    (4.20)& 4.08 (1.03) &--\\
	&\multirow{1}{*}{200}&2.89 (0.84)  &19.45 (3.92) & 2.82 (0.73) &--\\
	
	&\multirow{1}{*}{500}& 1.65 (0.41)&16.01 (3.57) & 1.64 (3.54)&--   \\
		&\multirow{1}{*}{1000}&1.17 (0.23) &  14.52 (3.51)  & 1.15 (0.21)&--\\
	\midrule 
	\multirow{6}{*}{L}& \multirow{1}{*}{20}  & 15.56 (3.66)& 51.48 (8.07) & 16.34 (4.06)&--\\

	&\multirow{1}{*}{50}& 10.79 (3.13)&45.57 (6.43)&10.69 (3.17)&--\\

	&\multirow{1}{*}{100}& 7.97 (2.07) &35.02 (5.79)& 7.41 (2.08) &--\\

	&\multirow{1}{*}{200}& 5.30 (1.41) &28.60 (5.21) &5.17 (1.43) &-- \\
	&\multirow{1}{*}{500}& 3.23 (0.76) &21.28 (5.03) & 3.05 (0.73) &--\\
	&\multirow{1}{*}{1000}& 1.39 (1.28) &  15.51 (4.80)& 1.96 (2.10) &--  \\

	\midrule
\multirow{8}{*}{M}&\multirow{1}{*}{20}&2.01 (2.97)& 5.29 (3.53)&--&2.62 (2.20) \\

&\multirow{1}{*}{50}&1.39 (1.28) & 3.42 (2.50) &--& 1.96 (2.10)   \\
&\multirow{1}{*}{100}&1.39 (0.91) &3.01 (3.02) &--&    1.42 (1.27)  \\

&\multirow{1}{*}{200}&0.82 (0.86) &  2.83 (3.19)&--&  1.11 (1.12)  \\
&\multirow{1}{*}{500}&0.82 (1.87) & 2.27 (3.82) &--& 0.95 (1.97)  \\
&\multirow{1}{*}{1000}& 0.69 (0.58)& 1.82 (1.59) &--& 0.81 (0.69)\\
&\multirow{1}{*}{2000}&0.43 (1.43) & 0.92 (3.23)&--& 0.39 (1.18)\\
&\multirow{1}{*}{5000}&0.25 (0.68) & 0.75 (1.80)& --&0.30 (1.25)\\
	\bottomrule
\end{tabular}}
		\label{ErrGamma1}
\end{table} 
Table \ref{ErrGamma1} presents the smoothing study in the multivariate setup ($d=3$ and 5) with respect to densities K, L and M.  Sample sizes of $n=20,50,100,200,500,1000$ are considered  for each model. In addition, sample sizes $n\in\{2000,5000\}$ are added for dimension $d = 5$ because $n = 1000$ is moderate or small in this setup.  Thus, for all sample sizes,  the ISE values show the superiority of the pure standard gamma kernel closely followed by the combined gamma version and a bit further the pure modified one.

\section{Illustrative applications}\label{sec:Illustration}

The proposed nonparametric method is here illustrated on four examples of univariate, bivariate and trivariate datasets. In order to smooth the joint distributions, the multivariate gamma kernel estimator \eqref{gammaestimator} and \eqref{mix_gammaestimator} with  Bayesian adaptive bandwidths are used. Unless otherwise specified, the retained parameters for the prior model on bandwidths are always $\alpha=n^{2/5}$ and $\beta_\ell=1$ for all $\ell=1,\ldots,d$. Once again, the boundary region in \eqref{ro} are selected in the range of $h_\ell\approx 10^{-8}$  with the modified  and combined gamma kernels \eqref{mix_gammaestimator}, and for all $\ell=1,\ldots,d$. 

 Besides the graphical method to compare estimators, we here use a numerical  procedure according to Filipone and Sanguinetti~\cite{FiliSAngui11} and \cite{SK20} for our four illustrations examples.  Precisely, we first sampled subset of size $m_n$ for each sample size $n$: $m_n \in 25, 50, 75, 100$ for the first univariate dataset with $n=119$, $m_n \in 20, 30, 40, 50$ for the second univariate with $n=66$,  $m_n \in 100, 150, 200, 250$ for the bivariate with $n=272$, and $m_n \in 10, 15, 20, 25$ for the trivariate with  $n=42$. Then,  we used the remaining data to compute the average log-likelihood. This experience is repeated 100 times for each $m_n$.
 
\begin{figure}[!ht]
	\mbox{
		
		\subfloat[(a1)]{	\resizebox*{6.7cm}{!}{\includegraphics{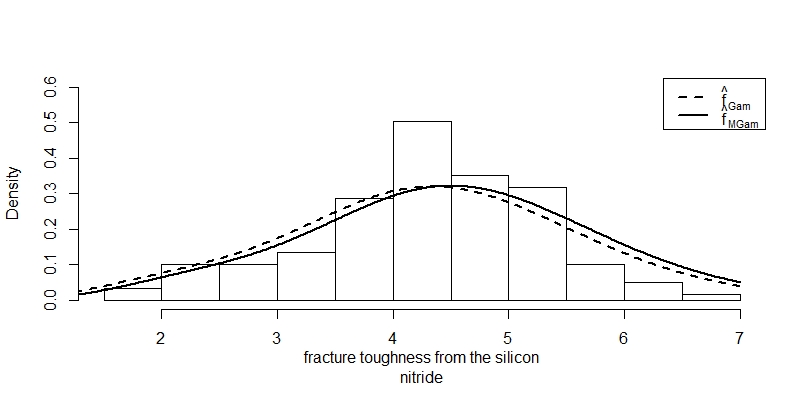}} }\hspace{1.5pt}
		\subfloat[(a2)]{\resizebox*{6.7cm}{!}{\includegraphics{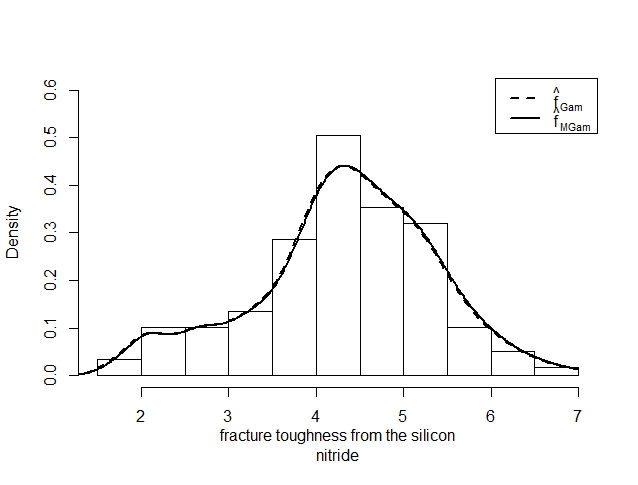}} }	
	}
	\mbox{
		
		\subfloat[(b1)]{	\resizebox*{6.7cm}{!}{\includegraphics{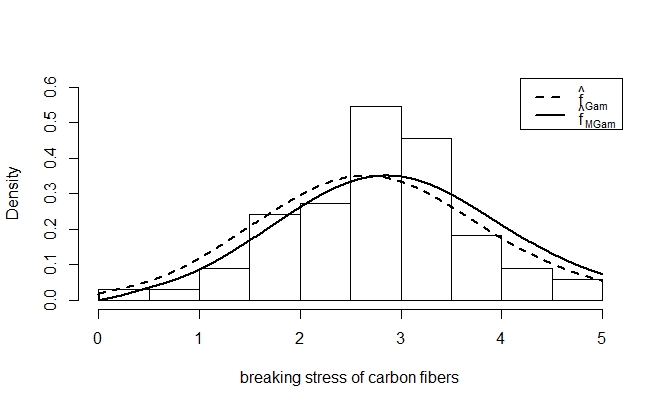}} }\hspace{1.5pt}
		\subfloat[(b2)]{\resizebox*{6.7cm}{!}{\includegraphics{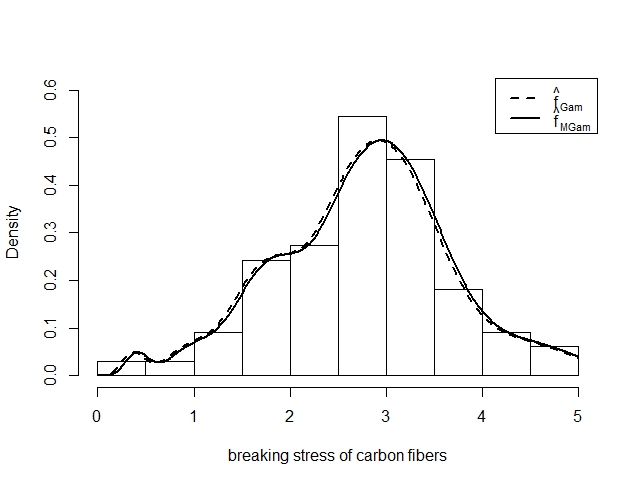}} }	
	}	
	
	\caption{Histograms    with their corresponding smoothings according to  Bayes selectors of bandwidths for univariate standard and modified gamma kernels: (a1--b1) $\alpha=n^{2/5}$ and $\beta_1=1$; (a2--b2) $\alpha=n^{4/5}$ and $\beta_1=1$.}
	\label{carbon_and_silicon}
\end{figure}

The first univariate example  concerns fracture toughness data (in the units of MPa $m^{1/2}$) of alumina ($Al_{2}O_{3}$). These  data have already been discussed in literature by \cite{NadarajahKotz2007}   and  \cite{ArshadElal2021}   and are available at  \href{http://www.ceramics.nist.gov/srd/summary/ftmain.htm}{http://www.ceramics.nist.gov/srd/summary/ftmain.htm}. Table \ref{stat_desc_univ}(a) provides a descriptive summary of data, with the standard deviation (SD) and skewness (CS), among others indicators. These nonnegative data are left-skewed and suggest the use of asymmetric kernels. Fig. \ref{carbon_and_silicon}(a1)–(a2) show the histogram and the smoothed distributions: the  solid line for the modified gamma kernel estimator and dashed line for standard one. We observe quasi-similar performances for both methods. The chosen parameters $\alpha=\alpha_n=n^{2/5}$ are not necessarily the optimal; see, e.g. \cite{SK20}. Thus,   Fig. \ref{carbon_and_silicon}(a2) obtained with $\alpha_n=n^{4/5}$ and $\beta_1=1$, provides a more precise smoothing of data. These previous statements are corroborated by the numerical results of the cross-validated log-likelihood method (Part (a) of Table \ref{log_univ}.)
\begin{table}[!htbp]
	\begin{center}
		\caption{Descriptive summary of the analyzed fracture toughness (a)  and tensile strength of carbon fibers (b) data.} \label{stat_desc_univ}
		\begin{tabular}{lrrrrrrr}
			\toprule
			Data set &\multicolumn{1}{c}{$n$}  &\multicolumn{1}{c}{max.}&\multicolumn{1}{c}{min.} &\multicolumn{1}{c}{median} &\multicolumn{1}{c}{mean}&\multicolumn{1}{c}{SD}&\multicolumn{1}{c}{CS}\\\hline
			Fracture roughness (a)	& 119 & 6.810  & 1.680 &4.380  &4.325  &   1.018&-0.411\\
			Breaking stress (b) & 66 & 4.900 & 0.390 &2.835  &2.760  &  
			0.891 &-0.128\\		
			\bottomrule
		\end{tabular}
	\end{center}
\end{table}
\begin{table}
	\caption{Mean average log-likelihood and their standard errors (in parentheses) for the analyzed  fracture toughness (a),  and tensile strength of carbon fibers (b) data, and  based on 100 replications.} 
	\begin{center}
		{		\begin{tabular}{rrrrrr}
				\toprule
				&& \multicolumn{2}{c}{$\alpha=n^{2/5}$, $\beta_1=1$} & \multicolumn{2}{c}{$\alpha=n^{4/5}$, $\beta_1=1$} \\
				\cmidrule(lr){3-4} \cmidrule(lr){5-6}
				
				&$m_n$&\multicolumn{1}{c}{Gamma}&\multicolumn{1}{c}{Modified gamma }& \multicolumn{1}{c}{Gamma}&\multicolumn{1}{c}{Modified gamma }   \\\midrule
				
				\multirow{4}{*}{(a)}	& \multirow{1}{*}{25}   &$-30.33$ (1.58)&$-28.82$ (0.75)&    $-29.36$ (2.67) & $-27.18$ (1.13)  \\	 
				&\multirow{1}{*}{50}&$-67.62$ (2.61) &$-73.71$ (1.36)&      $-64.66$ (3.76)  & $-73.22$ (2.54)         \\ 
				&\multirow{1}{*}{75} &$-103.63$ (2.79)&$-106.98$ (0.38) &    $-99.72$ (4.95)   &  $-102.74$ (1.22)   \\
				&\multirow{1}{*}{100}	&$-140.39$ (2.20)&$-139.75$ (0.49)&$-135.53$ (3.06) &$-132.59$ (0.78)   \\	
				
				\midrule
				
				\multirow{4}{*}{(b)}	& \multirow{1}{*}{20}	  &$-24.03$ (1.72)&$-22.57$ (0.76)&  $-22.67$ (4.14)&$-20.26$ (1.17)  \\	 
				&\multirow{1}{*}{30}&$-37.89$ (2.06)  & $-35.73$ (0.86)&$-36.97$ (3.63)&$-32.41$ (1.28)\\ 
				&\multirow{1}{*}{40} &$-51.55$ (2.05) &  $-49.28$ (0.70) &$-52.26$ (3.61)&$-45.52$ (1.03) \\
				&\multirow{1}{*}{50}&$-65.00$ (2.29)& $-66.38$ (0.46)&$-65.64$  (4.09)&$-63.51$ (2.66)\\	
				\bottomrule
				
		\end{tabular}}
		\label{log_univ}
	\end{center}
\end{table}

The second univariate example comes from  tensile strength of carbon fibers and has been discussed earlier in \cite{MichelePadgett2006}  and  \cite{ArshadElal2021}. Table \ref{stat_desc_univ}(b), Table \ref{log_univ}(b) and Fig.~\ref{carbon_and_silicon}(b1--b2) give successively the descriptive summary, the mean average log-likelihood and, finally, the histogram and the smoothed distributions of the dataset. Again, we observe the nonnegative nature, the unimodality and negative skewness of the distribution. Table \ref{log_univ}(b) shows the superiority of the modified gamma estimator with Bayesian adaptive bandwidths. Finally, Fig. \ref{carbon_and_silicon}(b2) and Table~\ref{log_univ}(b) obtained with $\alpha_n=n^{4/5}$ and $\beta_1=1$, give  better smoothed distributions and numerical results, respectively.

\begin{figure}[!htbp]
	\mbox{
		
		\subfloat[(Contour gamma)]{	\resizebox*{6.cm}{!}{\includegraphics{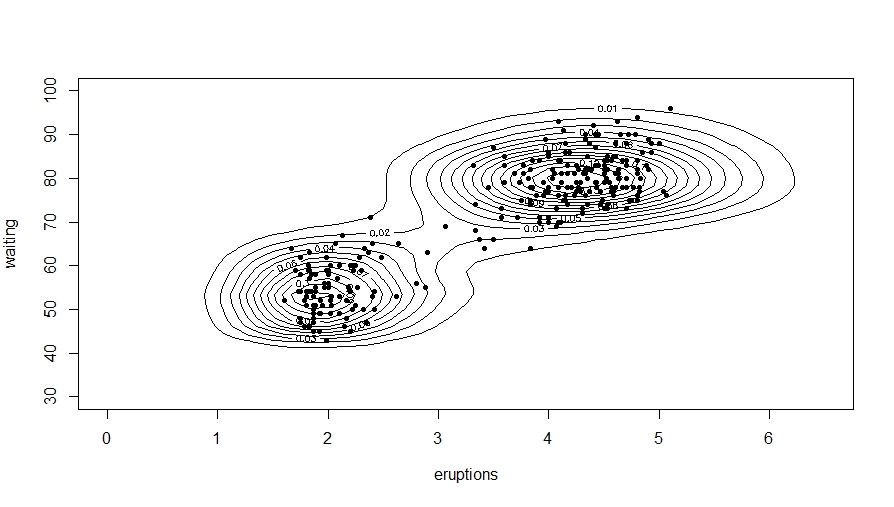}} }\hspace{2pt}
		\subfloat[(Surface gamma) ]{\resizebox*{8.cm}{!}{\includegraphics{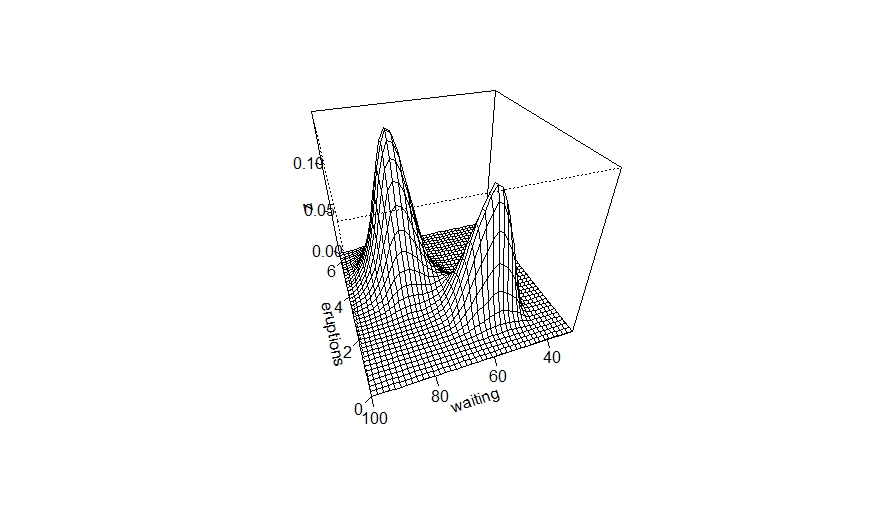}} }	
	}
	
	\mbox{
		
		\subfloat[(Contour-modified gamma)]{	\resizebox*{6.cm}{!}{\includegraphics{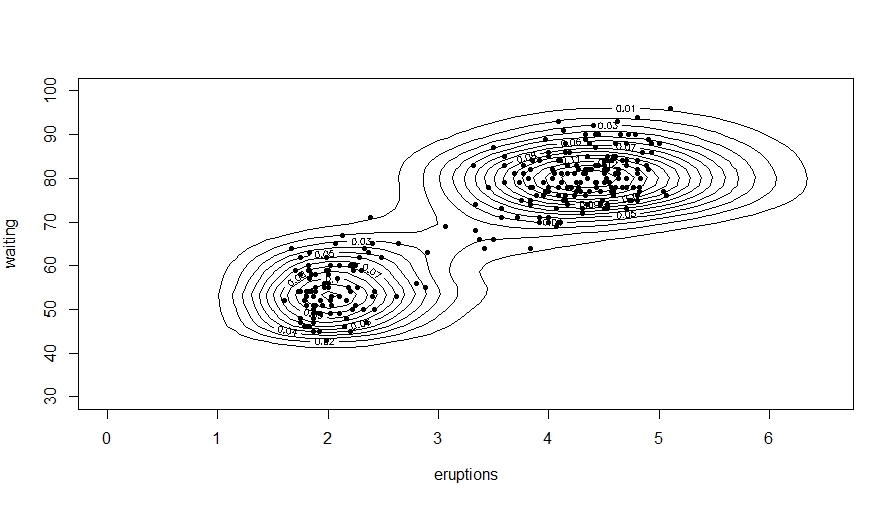}} }\hspace{2pt}
		\subfloat[(Surface-modified gamma) ]{\resizebox*{8.cm}{!}{\includegraphics{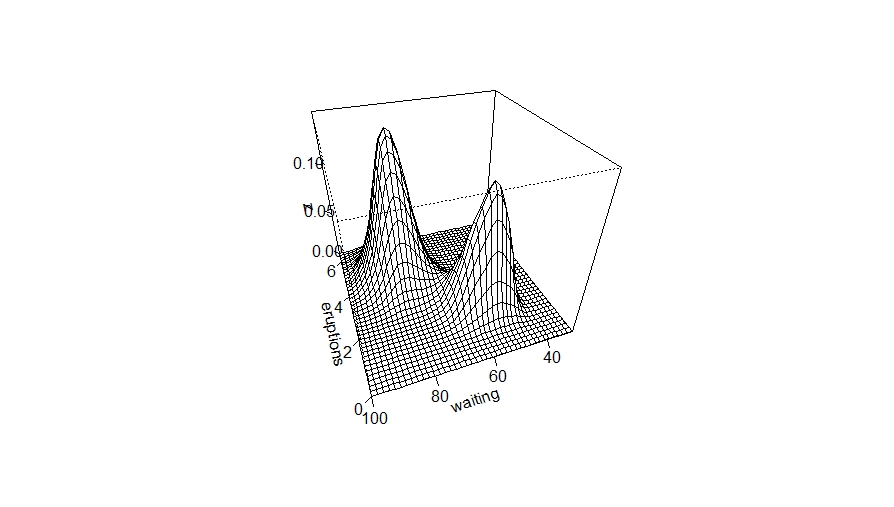}}} 
		
	}
	\mbox{ 
		
		\subfloat[(Contour-combined gamma)]{	\resizebox*{6.cm}{!}{\includegraphics{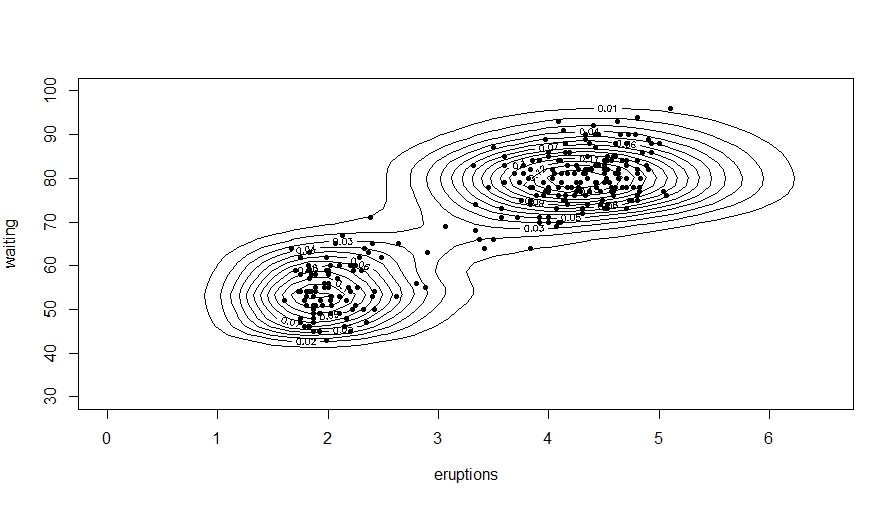}} }\hspace{2pt}
		\subfloat[(Surface-combined gamma)]{\resizebox*{8.cm}{!}{\includegraphics{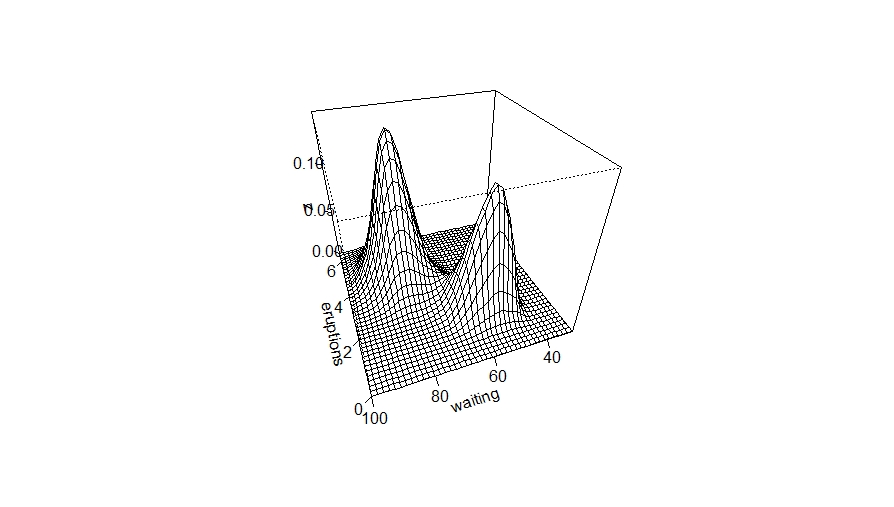}} }
	}
	
	\caption{Contour (left) and surface (right) plots of smoothed distribution from the Old Faithful geyser dataset according to Bayes selectors of bandwidths vector $\mathbf{h}$ with $\alpha=n^{2/5}$ and $\beta_1=\beta_2=1$.}
	\label{Old_faithful_data_est}
\end{figure}

The third bivariate illustration is from the Old Faithful geyser data already discussed in literature (see, e.g., \cite{AzzBow90}, \cite{FiliSAngui11} and recently \cite{SK20})  and available in the {\it Datasets Library} of the \textsf{R} software \cite{R20}. Data concern $n=272$ measurements of the eruption for the Old Faithful geyser in Yellowstone National Park, Wyoming, USA. The two covariates represent, in  minutes, the {\it waiting time} between eruptions and the {\it duration} of the eruption of these nonnegative real data.  Fig.~\ref{Old_faithful_data_est} show contour and surface plots of the smoothed distribution with Bayesian bandwidth selections for the Old Faithful geyser dataset. The points represent the scatter plots and the solid lines are the contour plot estimates. The smoothing of all methods are quite similar and point out  the bimodality of the Old Faithful geyser data. In contrast, the combined gamma method is the best  according to Table~\ref{log_Old_Geyser} of average log-likelihood.
\begin{table}
	\caption{Mean average log-likelihood and their standard errors (in parentheses) for  Old Faithful Geyser data based on 100 replications with $\alpha=n^{2/5}$ and $\beta_1=\beta_2=1$.} 
	\begin{center}
		{		\begin{tabular}{rrrrrr}
				\toprule
				&$m_n$&\multicolumn{1}{c}{Gamma}&\multicolumn{1}{c}{$Gam^{1}\times Mgam^{1}$ }&\multicolumn{1}{c}{Modified gamma }   \\\midrule
				& \multirow{1}{*}{100} &$-111.28$ (2.47)&$-110.23 $ (2.41)&$-135.80$ (6.18)   \\	 
				&\multirow{1}{*}{150}&$-329.50$ (54.69)  &$-327.24$ (54.46)& $-439.17$ (76.63) \\ 
				&\multirow{1}{*}{200}&$-570.01$ (5.42) &$-565.86$ (2.01)&  $-786.16$ (26.11)  \\
				&\multirow{1}{*}{250}&$-792.41$ (4.95) & $-792.71$ (1.26)& $-1130.75$ (33.16)	  \\ 		

				\bottomrule
		\end{tabular}}
		\label{log_Old_Geyser}
	\end{center}
\end{table}

The fourth trivariate illustration is done through the  drinking water pumps dataset which is recently used in~\cite{SK20, KS21} for nonparametric and semiparametric approaches, respectively. It concerns three measurements (with $n=42$) of drinking water pumps installed in the Sahel. 
The three variables refer successively to the {failure times} (in months) already used in Touré {\it et al.}~\cite{ToureEtal20}, the distance (in kilometers) between each water pump and the repair center, and finally the average volume  (in m$^3$) of water per~day. According to the average log-likelihood, the modified gamma kernel appears to be the best in general followed by the combined versions; see Table \ref{water_pump_est}.

\begin{table}
	\caption{Mean average log-likelihood and their standard errors (in parentheses) for drinking water pumps data based on 100 replications with $\alpha=n^{2/5}$ and $\beta_1=\beta_2=\beta_3=1$.} 
		\begin{center}
	{		\begin{tabular}{rrrrrr}
				\toprule
				&$m_n$&\multicolumn{1}{c}{Gamma}&\multicolumn{1}{c}{$Gam^{1}\times Mgam^{2}$ }&\multicolumn{1}{c}{$Gam^{2}\times Mgam^{1}$ }&\multicolumn{1}{c}{Modified gamma }   \\\midrule
				& \multirow{1}{*}{10}&$-190.59$ (17.02) & $-222.92$ (68.91)& $-222.96$ (69.16) &$-212.59$ (52.52)	 \\	 
				&\multirow{1}{*}{15}&$ -271.19$ (20.76) &$-269.25$ (35.09)& $-269.38$ (35.56) & $-264.66$ (23.66)  \\ 
				&\multirow{1}{*}{20}&$-344.26$ (17.60)  &$-320.64$ (29.11)&$-320.27$ (29.36) & $-320.25$ (25.92)\\
				&\multirow{1}{*}{25}  &$-426.54$ (18.67)&$-364.73$ (28.54)&$-364.65$ (28.69)  &$-363.62$ (23.95) \\ 		
		
				\bottomrule
	\end{tabular}}
	\label{water_pump_est}
		\end{center}
\end{table}

\section{Concluding remarks}
\label{sec:conclusion}

In this paper, we have proposed a semiparametric smoother method using combined gamma kernels and Bayesian adaptive selector of bandwidth vector. This multivariate kernel is the product of univariate standard and modified gamma kernels chosen according to each univariate margins of data. The new  method has also interesting asymptotic properties.  We have studied  the performance of the pure modified gamma kernel with a Bayesian adaptive bandwidths for nonnegative density smoothing. Efficient posterior and Bayes estimators of the bandwidths vector under quadratic loss function have been explicitly derived.

Simulation studies and four real datasets have highlighted the good performances of the proposed approaches (which are pure combined and pure modified) for the nonparametric estimation of the bandwidth in terms of ISE and log-likelihood criteria. As expected, the multiple combined gamma kernel smoothers are interesting compromises between the complete standard and modified gamma ones. It is even the best in some cases. For an efficient practical use, we recommend to check the shape of the univariate margins of the multivariate data in order to make the best choice (combined, standard or modified gamma kernel) accordingly. Future work in progress is devoted to a Bayesian local bandwidths for discrete multivariate associated kernel for semiparametric smoothers; see, e.g.,  \cite{KS21} in multivariate continuous case.

\section*{Appendix: Proofs of propositions}\label{sec:Appendix}

\textbf{Proof of Proposition \ref{PropBiasVarf(x,0)}.}
	Since one has $$\mathrm{Bias}[\widehat{f}_n(\boldsymbol{x})]=p_{d}(\boldsymbol{x};\widehat{\boldsymbol{\theta}}_n)\mathbb{E}[\widehat{w}_n(\boldsymbol{x})]-f(\boldsymbol{x}) \mbox{ and } \mathrm{var}[\widehat{f}_n(\boldsymbol{x})]=[p_{d}(\boldsymbol{x};\widehat{\boldsymbol{\theta}}_n)]^2 \mathrm{var} [\widehat{w}_n(\boldsymbol{x})],$$ it is enough to calculate $\mathbb{E}[\widehat{w}_n(\boldsymbol{x})]$ and $\mathrm{var}[\widehat{w}_n(\boldsymbol{x})]$ using $\widehat{w}_n(\boldsymbol{x})=n^{-1}\sum_{i=1}^n \mathbf{G}_{\boldsymbol{x},\mathbf{h},\ell}(\mathbf{X}_i)/p_{d}(\mathbf{X}_i;\widehat{\boldsymbol{\theta}}_n)$ for all $\boldsymbol{x}\in\mathbb{T}_d^+$ and $\mathbf{G}_{\boldsymbol{x},\mathbf{h},\ell}=\left(\prod_{s=1}^{d-\ell}G_{x_{s},h_{s}}\right)\left(\prod_{r=1}^{\ell}G_{\rho(x_{r};h_{r}),h_{r}}\right)$ from \eqref{parti_cases}. Indeed, one successively has
	\begin{align*}
	\mathbb{E} \left[\widehat w_n(\boldsymbol{x})\right]
	&= \mathbb{E}\left[\mathbf{G}_{\boldsymbol{x},\mathbf{h},\ell}(\mathbf{X}_{1})/p_{d}(\mathbf{X}_1;\widehat{\boldsymbol{\theta}}_n)\right]\\
	&=\int_{[0,\infty)^d}\mathbf{G}_{\boldsymbol{x},\mathbf{h},\ell}(\mathbf{u})\left[p_{d}(\mathbf{u};\widehat{\boldsymbol{\theta}}_n)\right]^{-1}f(\mathbf{u})d\mathbf{u}
	=\mathbb{E}\left[w\left({\cal G}_{\boldsymbol{x},\mathbf{h},\ell}\right)\right]\\
	&= w(\boldsymbol{x})+\left[\displaystyle \sum_{r=1}^{d-\ell}h_r \frac{\partial w}{\partial x_r}(\boldsymbol{x})+\displaystyle \sum_{r=1}^{d-\ell}\frac{1}{2}\left(x_r h_r+2h_r^2 \right)\frac{\partial^2 w}{\partial x_r^2}(\boldsymbol{x})+\frac{1}{2}\displaystyle \sum_{s=1}^{\ell}x_s h_s \frac{\partial^2 w}{\partial x_s^2}(\boldsymbol{x})\right] 
	 +\left(1+o\left\{\displaystyle \sum_{j=1}^{d}h_j^2 \right\}\right),
	\end{align*}
	which leads to the  result of $\mathrm{Bias}[\widehat{f}_n(\boldsymbol{x})]$.
	 
	About the variance term, $f$ being bounded leads to $\mathbb{E}\left[\mathbf{G}_{\boldsymbol{x},\mathbf{h},\ell}(\mathbf{X}_{j})\right]=O(1)$.  Also, we denote by $\nabla f(\boldsymbol{x})$ and $\mathcal{H}f(\boldsymbol{x})$ the gradient vector and the corresponding Hessian matrix of the function $f$ at $\boldsymbol{x}$, respectively. It successively follows:  
	\begin{align*}
	\mathrm{var} \left[\widehat w_n(\boldsymbol{x})\right]
	&= \frac{1}{n}\mathrm{var}\left[\mathbf{G}_{\boldsymbol{x},\mathbf{h},\ell}(\mathbf{X}_{1})/p_{d}(\mathbf{X}_1;\widehat{\boldsymbol{\theta}}_n)\right]\\
	&=\frac{1}{n}\left[\int_{[0,\infty)^d}\mathbf{G}^2_{\boldsymbol{x},\mathbf{h},\ell}(\mathbf{u})[p_{d}(\mathbf{u};\widehat{\boldsymbol{\theta}}_n)]^{-2}f(\mathbf{u})d\mathbf{u}+O(1)\right]\\
	&=\frac{1}{n} \int_{[0,\infty)^d}\mathbf{G}^2_{\boldsymbol{x},\mathbf{h},\ell}(\mathbf{u})[p_{d}(\mathbf{u};\widehat{\boldsymbol{\theta}}_n)]^{-2}\begin{pmatrix}f(\boldsymbol{x})+(\boldsymbol{x}-\mathbf{u})^T\nabla f(\boldsymbol{x})\\+\frac{1}{2}(\boldsymbol{x}-\mathbf{u})^T{\cal H}f(\boldsymbol{x}) (\boldsymbol{x}-\mathbf{u})\\
	+o\left[\left(||\boldsymbol{x}-\mathbf{u}||^2\right)\right]
	\end{pmatrix}
	d\mathbf{u}\\
	&=\frac{1}{n}f(\boldsymbol{x})[p_{d}(\boldsymbol{x};\widehat{\boldsymbol{\theta}}_n)]^{-2}||\mathbf{G}_{\boldsymbol{x},\mathbf{h},\ell}||_2^2++o\left(n^{-1}\displaystyle\prod_{j=1}^{d}h_j^{-1/2}\right)\\
	&=\frac{1}{n}f(\boldsymbol{x})[p_{d}(\boldsymbol{x};\widehat{\boldsymbol{\theta}}_n)]^{-2}\displaystyle \prod_{k=1}^{d-\ell}\left(\frac{\Gamma(1+2x_k/h_k)}{2^{1+2x_k/h_k}\Gamma(1+x_k/h_k)}h_{k}^{-1}\right)\prod_{s \in \mathbb{I}^{}_{2}}\left(\frac{\Gamma(1+\lambda_{s}^2/2)}{2^{1+\lambda_{s}^2/2}\Gamma(1+\lambda_{s}^2/4)}h_{s}^{-1}\right)\\
	&\quad \prod_{j \in \mathbb{I}_{2}^c}\left(\frac{1}{2\pi^{1/2}}h_j^{-1/2}x_j^{-1/2}\right)+o\left(n^{-1}\displaystyle\prod_{j=1}^{d}h_j^{-1/2}\right)\\
	&=\frac{1}{n}f(\boldsymbol{x})[p_{d}(\boldsymbol{x};\widehat{\boldsymbol{\theta}}_n)]^{-2}\displaystyle \prod_{k \in \mathbb{I}^{}_{1}}\left(\frac{\Gamma(1+2\lambda_{k})}{2^{1+2\lambda_{k}}\Gamma(1+\lambda_{k})}h_{k}^{-1}\right)\prod_{s \in \mathbb{I}^{}_{2}}\left(\frac{\Gamma(1+\lambda_{s}^2/2)}{2^{1+\lambda_{s}^2/2}\Gamma(1+\lambda_{s}^2/4)}h_{s}^{-1}\right)\nonumber\\
	&\quad \prod_{j \in \mathbb{I}_{}^c}\left(\frac{1}{2\pi^{1/2}}h_j^{-1/2}x_j^{-1/2}\right)+o\left(n^{-1}\displaystyle\prod_{j=1}^{d}h_j^{-1/2}\right),
	\end{align*}  and the desired result of $\mathrm{var}[\widehat{f}_n(\boldsymbol{x})]$ is therefore deduced. $\blacksquare$

\textbf{Proof of Proposition \ref{theo1}.} 	(i) Let us represent  $\pi(\mathbf{h}_{i}\mid\mathbf{X}_{i})$ of (\ref{bayesrule}) as the ratio of  $N(\mathbf{h}_{i}\mid\mathbf{X}_{i}):=\widehat{f}_{n,\mathbf{h}_i,-i}(\mathbf{X}_i) \pi (\mathbf{h}_{i})$ and  $\int_{[0, \infty)^{d}}N(\mathbf{h}_{i}\mid\mathbf{X}_{i})d\mathbf{h}_{i}$.
	From $(\ref{equ5})$ and $(\ref{prior})$ the numerator is first equal to 
	\begin{align}
	N(\mathbf{h}_{i}\mid\mathbf{X}_{i})&= \left(\frac{p_{d}(\mathbf{X}_i;\widehat{\boldsymbol{\theta}}_n)}{n-1}\sum_{j=1,j\neq i}^{n}
	\frac{1}{p_{d}(\mathbf{X}_{j};\widehat{\boldsymbol{\theta}}_n)} \prod_{\ell=1}^{d} G_{\rho(X_{i\ell},h_{i\ell}),h_{i\ell}}(X_{j\ell})\right)\left(\prod_{\ell=1}^{d} \frac{\beta_{\ell}^{\alpha}}{\Gamma(\alpha)}h_{i\ell}^{-\alpha-1}\exp(-\beta_{\ell}/h_{i\ell})\right)\nonumber\\
	&=\frac{p_{d}(\mathbf{X}_i;\widehat{\boldsymbol{\theta}}_n)}{(n-1)[\Gamma(\alpha)]^{d}}\sum_{j=1,j\neq i}^{n}
	\frac{1}{p_{d}(\mathbf{X}_{j};\widehat{\boldsymbol{\theta}}_n)} \prod_{\ell=1}^{d} \frac{G_{\rho(X_{i\ell},h_{i\ell}),h_{i\ell}}(X_{j\ell})}{\beta_{\ell}^{-\alpha}h_{i\ell}^{\alpha+1}}\exp(-\beta_{\ell}/h_{i\ell}).\label{equ8}
	\end{align}	
	
	From $(\ref{ro})$, consider the following partition $\mathbb{I}_{\mathbf{X}_i}$ and $\mathbb{I}_{\mathbf{X}_i}^{c}$ of $\{1,2,...,d\}$. For $X_{ik} \in [0,2h_{ik})$ with $k\in \mathbb{I}_{\mathbf{X}_i}$, the function $n\mapsto \left[X_{ik}/2h_{ik}(n)\right]^{2}$ is bounded and then there exists a constant $\lambda_{ik}>0$ such that $(X_{ik}/2h_{ik})^{2} \rightarrow \lambda_{ik}$ as $n\rightarrow \infty$; see Chen~\cite[pp. 474-475]{C00}. Using successively ($\ref{gam2}$) and ($\ref{ro}$) with the behavior $(X_{ik}/2h_{ik})^{2}\simeq \lambda_{ik}$ as $n\rightarrow \infty$, the term of product on $\mathbb{I}_{\mathbf{X}_i}$ in $(\ref{equ8})$ can be expressed as follows
	\begin{align}\label{N1}
	\frac{G_{\rho(X_{ik};h_{ik}),h_{ik}}(X_{jk})}{\beta_{k}^{-\alpha}h_{ik}^{\alpha+1}}\exp(-\beta_{k}/h_{ik})&=\frac{X_{jk}^{(X_{ik}/2h_{ik})^{2}}\exp(-X_{jk}/h_{ik})}{h_{ik}^{1+(X_{ik}/2h_{ik})^{2}}
		\Gamma[1+(X_{ik}/2h_{ik})^{2}]\beta_{k}^{-\alpha}h_{ik}^{\alpha+1}}\exp(-\beta_{k}/h_{ik})\nonumber\\
	&\simeq\frac{X_{jk}^{\lambda_{ik}}\exp[-(X_{jk}+
		\beta_{k})/h_{ik}]}{h_{ik}^{\lambda_{ik}+\alpha+2}\beta_{k}^{-\alpha}\Gamma(1+\lambda_{ik})}\nonumber\\
	&=\frac{ \Gamma (\lambda_{ik}+ \alpha +1) X_{jk}^{\lambda_{ik}}}{\beta_{k}^{-\alpha}\Gamma (\lambda_{ik}+1)(X_{jk}+\beta_{k})^{\lambda_{ik}+\alpha+1}} \times \frac{(X_{jk}+\beta_{k})^{\lambda_{ik}+\alpha+1}\exp[-(X_{jk}+
		\beta_{k})/h_{ik}]}{h_{ik}^{\lambda_{ik}+\alpha+2}\Gamma(\lambda_{ik}+\alpha+1)}\nonumber\\
	&=A_{ijk}\,IG_{\lambda_{ik}+\alpha+1,X_{jk}+
		\beta_{k}}(h_{ik}),
	\end{align}
	with $A_{ijk}(\alpha,\beta_k)= [ \Gamma (\lambda_{ik}+ \alpha +1) X_{jk}^{\lambda_{ik}}]/[\beta_{k}^{-\alpha}\Gamma (\lambda_{ik}+1)(X_{jk}+\beta_{k})^{\lambda_{ik}+\alpha+1}]$ and  $IG_{\lambda_{ik}+\alpha+1,X_{jk}+\beta_{k}}(h_{ik})$ comes from $(\ref{prior})$.
	
	Consider the largest part $\mathbb{I}^{c}_{\mathbf{X}_i}=\left \{\ell \in \{1,\ldots,d\}~;X_{i\ell} \in [2h_{i\ell},  \infty)\right \}$. Following again \citet[pp. 474-475]{C00}, we assume that for all $X_{i\ell}\in [2h_{i\ell}, \infty)$  one has $X_{i\ell}/h_{i\ell} \rightarrow \infty$ as $n\rightarrow \infty$ for all $\ell \in \{1,2,\ldots,d\}$. From  \eqref{gam2}, \eqref{ro}, the Sterling formula $\Gamma(z+1)\simeq\sqrt{2\pi}z^{z+1/2}\exp(-z)$ as $z\rightarrow \infty$,  and the well-known property $\Gamma(z)=z^{-1}\Gamma(z+1)$ for $z>0$, the term  $(\ref{equ8})$ can be successively calculated as
	\begin{align}\label{N2}	\frac{G_{\rho(X_{i\ell};h_{i\ell}),h_{i\ell}}(X_{j\ell})}{\beta_{\ell}^{-\alpha}h_{i\ell}^{\alpha+1}}\exp(-\beta_{\ell}/h_{i\ell})&=\frac{X_{j\ell}^{(X_{i\ell}/h_{i\ell})-1}\exp(-X_{j\ell}/h_{i\ell})}{h_{i\ell}^{X_{i\ell}/h_{i\ell}}
		\Gamma(X_{i\ell}/h_{i\ell})\beta_{\ell}^{-\alpha}h_{i\ell}^{\alpha+1}}\exp(-\beta_{\ell}/h_{i\ell})\nonumber\\
	&=\frac{X_{j\ell}^{(X_{i\ell}/h_{i\ell})-1}\exp(-X_{j\ell}/h_{i\ell})}{h_{i\ell}^{X_{i\ell}/h_{i\ell}}
		(X_{i\ell}/h_{i\ell})^{-1}\Gamma(1+X_{i\ell}/h_{i\ell})\beta_{\ell}^{-\alpha}h_{i\ell}^{\alpha+1}}\exp(-\beta_{\ell}/h_{i\ell})\nonumber\\
	&=\frac{X_{j\ell}^{-1}}{\beta_{\ell}^{-\alpha}X_{i\ell}^{-1}}\frac{\exp[-(X_{j\ell}+
		\beta_{\ell}-X_{i\ell}\log X_{j\ell})/h_{i\ell}]}{h_{i\ell}^{(X_{i\ell}/h_{i\ell})+\alpha+2}\sqrt{2\pi}(X_{i\ell}/h_{i\ell})^{(X_{i\ell}/h_{i\ell})+1/2}\exp(-X_{i\ell}/h_{i\ell})}\nonumber\\
	&=\frac{X_{j\ell}^{-1}\Gamma(\alpha+1/2)}{\beta_{\ell}^{-\alpha}X_{i\ell}^{-1/2}\sqrt{2\pi}[C_{ij\ell}(\beta_\ell)]^{\alpha+1/2}} \times \frac{[C_{ij\ell}(\beta_\ell)]^{\alpha+1/2}\exp(-C_{ij\ell}(\beta_\ell)/h_{i\ell})}{h_{i\ell}^{\alpha+3/2}\Gamma(\alpha+1/2)}\nonumber\\
	&=B_{ij\ell}(\alpha,\beta_\ell)\,IG_{\alpha+1/2,C_{ij\ell}(\beta_\ell)}(h_{i\ell}),
	\end{align}
	with $B_{ij\ell}(\alpha,\beta_\ell)= [X_{j\ell}^{-1}\Gamma(\alpha +1/2)]/(\beta_{\ell}^{-\alpha}X_{i\ell}^{-1/2}\sqrt{2\pi}[C_{ij\ell}(\beta_\ell)]^{\alpha +1/2})$, $C_{ij\ell}(\beta_\ell)= X_{i\ell}\log (X_{i\ell}/X_{j\ell})+X_{j\ell}-X_{i\ell}+\beta_{\ell}$ and $IG_{\alpha+1/2,C_{ij\ell}(\beta_\ell)}(h_{i\ell})$ is given in $(\ref{prior})$.
	
	Combining  ($\ref{N1}$) and ($\ref{N2}$), the expression  $N(\textbf{h}_{i}\mid\textbf{X}_{i})$ in $(\ref{equ8})$ becomes
	\begin{align}\label{N}
	\hspace{-0.75cm}N(\textbf{h}_{i}\mid\textbf{X}_{i})&=\frac{p_{d}(\mathbf{X}_i;\widehat{\boldsymbol{\theta}}_n)}{(n-1)[\Gamma(\alpha)]^{d}}\sum_{j=1,j\neq i}^{n}	\frac{1}{p_{d}(\mathbf{X}_{j};\widehat{\boldsymbol{\theta}}_n)} \left(\prod_{k \in \mathbb{I}_{\mathbf{X}_i}}A_{ijk}(\alpha,\beta_k)\,IG_{\lambda_{ik}+\alpha+1,X_{jk}+\beta_{k}}(h_{ik})\right) \left(\prod_{\ell \in \mathbb{I}^{c}_{\mathbf{X}_i}}B_{ij\ell}(\alpha,\beta_\ell)\,IG_{\alpha+1/2,S_{ij\ell}}(h_{i\ell})\right). 
	\end{align}
	From $(\ref{N})$, the denominator is successively computed as follows
	\begin{align}\label{D}
	\int_{[0, \infty)^{d}} N(\textbf{h}_{i}\mid\textbf{X}_{i})\,d\textbf{h}_{i}\nonumber&=\frac{p_{d}(\mathbf{X}_i;\widehat{\boldsymbol{\theta}}_n)}{(n-1)[\Gamma(\alpha)]^{d}}\sum_{j=1,j\neq i}^{n}
	\frac{1}{p_{d}(\mathbf{X}_{j};\widehat{\boldsymbol{\theta}}_n)} \left(\prod_{k \in \mathbb{I}_{\mathbf{X}_i}}A_{ijk}(\alpha,\beta_k)\int_{0}^{\infty}IG_{\lambda_{ik}+\alpha+1,X_{jk}+
		\beta_{k}}(h_{ik})\,dh_{ik}\right)\nonumber\\
	&\quad\times \left(\prod_{\ell \in \mathbb{I}^{c}_{\mathbf{X}_i}}B_{ij\ell}(\alpha,\beta_\ell)\int_{0}^{\infty}IG_{\alpha+1/2,C_{ij\ell}(\beta_\ell)}(h_{i\ell})\,dh_{i\ell}\right)\nonumber\\
	&=\frac{p_{d}(\mathbf{X}_i;\widehat{\boldsymbol{\theta}}_n)}{(n-1)[\Gamma(\alpha)]^{d}}\sum_{j=1,j\neq i}^{n}
	\frac{1}{p_{d}(\mathbf{X}_{j};\widehat{\boldsymbol{\theta}}_n)} \left(\prod_{k \in \mathbb{I}_{\mathbf{X}_i}}A_{ijk}(\alpha,\beta_k)\right)\left(\prod_{\ell \in \mathbb{I}^{c}_{\mathbf{X}_i}}B_{ij\ell}(\alpha,\beta_\ell)\right)\nonumber\\
	&=\frac{p_{d}(\mathbf{X}_i;\widehat{\boldsymbol{\theta}}_n)}{(n-1)[\Gamma(\alpha)]^{d}}D_{i}(\alpha,\boldsymbol{\beta}),
	\end{align}	with $D_{i}(\alpha,\boldsymbol{\beta})=p_{d}(\mathbf{X}_i;\widehat{\boldsymbol{\theta}}_n)\sum_{j=1,j\neq i}^{n}\left(p_{d}(\mathbf{X}_{j};\widehat{\boldsymbol{\theta}}_n)\right)^{-1}\left(\prod_{k \in \mathbb{I}_{\mathbf{X}_i}}A_{ijk}(\alpha,\beta_k)\right)\left(\prod_{\ell \in \mathbb{I}^{c}_{\mathbf{X}_i}} B_{ij\ell}(\alpha,\beta_\ell)\right)$. Finally, the ratio of $(\ref{N})$ and $(\ref{D})$ leads to Part (i).
	
	(ii) We remind that the mean of the inverse gamma distribution  $\mathcal{IG}(\alpha,\beta_\ell)$ is  $\beta_\ell/(\alpha-1)$ and $\mathbb{E}(h_{i\ell}\mid\mathbf{X}_{i})=\int_{0}^{\infty} h_{i\ell}\pi(h_{im}\mid\mathbf{X}_{i})\,dh_{im}$ with $\pi(h_{im}\mid\mathbf{X}_{i})$ is the marginal distribution $h_{im}$ obtained by integration of $\pi(\mathbf{h}_{i}\mid\mathbf{X}_{i})$  for all components of  $\mathbf{h}_{i}$ except $h_{im}$. Then, $\pi(h_{im}\mid\mathbf{X}_{i})=\int_{[0, \infty)^{d-1}}\pi(\mathbf{h}_{i}\mid\mathbf{X}_{i})\,d\mathbf{h}_{i(-m)}$ where $d\mathbf{h}_{i(-m)}$ is the vector  $d\mathbf{h}_{i}$ without the $m^{th}$ component.
	If $m\in \mathbb{I}_{\mathbf{X}_i}$, one has
	\begin{equation*}
	\pi(h_{im}\mid\mathbf{X}_{i})=\frac{p_{d}(\mathbf{X}_i;\widehat{\boldsymbol{\theta}}_n)}{D_{i}(\alpha,\boldsymbol{\beta})}\sum_{j=1,j\neq i}^{n}
	\frac{1}{p_{d}(\mathbf{X}_{j};\widehat{\boldsymbol{\theta}}_n)}\left(\prod_{k \in \mathbb{I}_{\mathbf{X}_i}}A_{ijk}(\alpha,\beta_k)\right)\left(\prod_{\ell \in \mathbb{I}^{c}_{\mathbf{X}_i}}B_{ij\ell}(\alpha,\beta_\ell)\right)IG_{\alpha+1, X_{jm}+\beta_{m}}(h_{im})
	\end{equation*}
	and
	\begin{equation}\label{h1}
	\widehat{h}_{im}=\mathbb{E}(h_{im}\mid\mathbf{X}_{i})=\frac{p_{d}(\mathbf{X}_i;\widehat{\boldsymbol{\theta}}_n)}{D_{i}(\alpha,\boldsymbol{\beta})}\sum_{j=1,j\neq i}^{n}
	\frac{1}{p_{d}(\mathbf{X}_{j};\widehat{\boldsymbol{\theta}}_n)}\left(\prod_{k \in \mathbb{I}_{\mathbf{X}_i}}A_{ijk}(\alpha,\beta_k)\right)\left(\prod_{\ell \in \mathbb{I}^{c}_{\mathbf{X}_i}}B_{ij\ell}(\alpha,\beta_\ell)\right)\left(\frac{X_{jm}+\beta_{m}}{\lambda_{im}+\alpha}\right).
	\end{equation}
	If $m\in \mathbb{I}_{\mathbf{X}_i}^{c}$ and $\alpha>1/2$, one gets
	\begin{equation*}
	\pi(h_{im}\mid\mathbf{X}_{i})=\frac{p_{d}(\mathbf{X}_i;\widehat{\boldsymbol{\theta}}_n)}{D_{i}(\alpha,\boldsymbol{\beta})}\sum_{j=1,j\neq i}^{n}
	\frac{1}{p_{d}(\mathbf{X}_{j};\widehat{\boldsymbol{\theta}}_n)}\left(\prod_{k \in \mathbb{I}_{\mathbf{X}_i}}A_{ijk}(\alpha,\beta_k)\right)\left(\prod_{\ell \in \mathbb{I}^{c}_{\mathbf{X}_i}}B_{ij\ell}(\alpha,\beta_\ell)\right)IG_{\alpha+1/2,C_{ijm}(\beta_m)}(h_{im})
	\end{equation*}
	and
	\begin{equation}\label{h2}
	\widehat{h}_{im}=\mathbb{E}(h_{im}\mid\mathbf{X}_{i})=\frac{p_{d}(\mathbf{X}_i;\widehat{\boldsymbol{\theta}}_n)}{D_{i}(\alpha,\boldsymbol{\beta})}\sum_{j=1,j\neq i}^{n}
	\frac{1}{p_{d}(\mathbf{X}_{j};\widehat{\boldsymbol{\theta}}_n)}\left(\prod_{k \in \mathbb{I}_{\mathbf{X}_i}}A_{ijk}(\alpha,\beta_k)\right)\left(\prod_{\ell \in \mathbb{I}^{c}_{\mathbf{X}_i}}B_{ij\ell}(\alpha,\beta_\ell)\right)\left(\frac{C_{ijm}(\beta_m)}{\alpha-1/2}\right).
	\end{equation}
	Combining $(\ref{h1})$ and $(\ref{h2})$, we therefore get the closed expression of Part (ii). $\blacksquare$


\bibliographystyle{elsarticle-harv}

\section*{Supplementary material}
The additional material  contains   surface and contour plots of the  target densities E, F, G and H, and their corresponding smoothed densities from Section \ref{SSection_Biv} of the bivariate simulation studies.

\end{document}